\documentclass[11pt]{article}
\usepackage{titlesec,fullpage, lmodern}
\usepackage[english]{babel}
\usepackage[utf8]{inputenc}
\usepackage[T1]{fontenc}
\usepackage{graphicx, enumitem, xifthen}
\usepackage{mathrsfs,amsthm,amssymb,amsfonts,amsmath}
\usepackage[all]{xy} 

\usepackage[top=2cm,bottom=3cm,right=2.5cm,left=2.5cm]{geometry} 

\titleformat{\section}{\Large\sc}{\thesection.}{1em}{} 

\newtheoremstyle{plain} {\topsep}{\topsep}{\itshape}{}{\bf}{ $-$}{ }
{\thmname{#1}\thmnumber{ #2}\thmnote{ \normalfont(#3)}}

\newtheoremstyle{defin}
{\topsep}{\topsep}{\rm}{}{\bf}{ $-$}{ }
{\thmname{#1}\thmnumber{ #2}\thmnote{ \normalfont(#3)}}

\newtheorem{theoreme}{Theorem}[section]
\newtheorem{lemme}[theoreme]{Lemma}
\theoremstyle{defin}
\newtheorem{definition}[theoreme]{Definition}
\newtheorem{definitions}[theoreme]{Definitions}

\newtheorem{notation}[theoreme]{Notation}
\newtheorem{rmq}[theoreme]{Remark}
\newtheorem{rmqs}[theoreme]{Remarks}

\newenvironment{exemples}{ \refstepcounter{theoreme} \noindent
  \textbf{Examples \thetheoreme.} \rm}{}

\newenvironment{preuve}[1][]
{\pushQED{\qed}
\noindent \textsc{Proof\ifthenelse{\isempty{#1}}{}{ of Theorem #1}. }}
{ \unskip \nobreak \hfill \popQED\vspace{\itemsep}  }

\renewcommand{\ker}{\textrm{Ker}} \newcommand{\im}{\textrm{Im}}
\newcommand{\dom}{\textrm{Dom}} \renewcommand{\sp}{\textrm{sp}}
\newcommand{\vect}{\textrm{span}}

\newcommand{\Ast}[1]{\underset{#1}{\scalebox{2}{\raisebox{-0.2ex}{$\ast$}}}}%
\newcommand{\bigobot}[2][]{\overset{#1}{\underset{#2}{\scalebox{1.4}{$\bigcirc\kern-12.8pt\perp\,$}}}}

\newcommand{\NN}{\mathbb N}
\newcommand{\ZZ}{\mathbb Z}

\newcommand{\CC}{\mathbb C}

\newcommand{\cA}{\mathcal{A}} 
\newcommand{\cC}{\mathcal{C}} 
\newcommand{\cE}{\mathcal{E}} \newcommand{\cF}{\mathcal{F}}
 \newcommand{\cH}{\mathcal{H}}
\newcommand{\cI}{\mathcal{I}} \newcommand{\cJ}{\mathcal{J}}
 \newcommand{\cL}{\mathcal{L}}
\newcommand{\cM}{\mathcal{M}} 
 
\newcommand{\cQ}{\mathcal{Q}} 
 
\newcommand{\cU}{\mathcal{U}}

\newcommand{\0}{_{(0)}}
\renewcommand{\1}{_{(1)}}

\renewcommand{\d}{\mathrm{d}}  
 \newcommand{\I}{\mathscr I}
\newcommand{\cat}{\cC (E, \tau, (V_i)_{i \in \cI}, J, \xi_0) }
\renewcommand{\dot}{\,\cdot\,}
\renewcommand{\l}{\langle}
\renewcommand{\r}{\rangle}
\newcommand{\cf}{\textit{cf.~}}
\newcommand{\ie}{i.e. }
\newcommand{\tq}{\;\, ; \;}
\newcommand{\appl}[5][]{\begin{array}[t]{crcl}
{#1} &\!\! {#2} & \!\!\rightarrow\!\! & {#3} \!\!\\ 
&\!\! {#4} & \!\!\mapsto\!\! & {#5} \!\!
\end{array}}

\AtBeginDocument{} 
\AtBeginDocument{\setenumerate{label=(\alph*),leftmargin=1cm}} 
\makeatletter\@beginparpenalty=0 \makeatother 
\title{Quantum symmetry groups of Hilbert modules equipped with orthogonal filtrations}
\date{}
\author{Manon Thibault de Chanvalon}

\begin{document}
\maketitle
\begin{center}
\textit{{Laboratoire de~Math\'ematiques (UMR 6620)}, Universit\'e Blaise Pascal, \\ Complexe universitaire des C\'ezeaux, 63171 Aubi\`ere Cedex, France.}
\end{center}
\begin{center}
{\tt manon.thibault@math.univ-bpclermont.fr}
\vspace{1em}
\end{center}

\begin{abstract}
We define and show the existence of the quantum symmetry group of a Hilbert module equipped
with an orthogonal filtration. Our construction unifies the constructions of Banica-Skalski's quantum symmetry group
of a $C^*$-algebra equipped with an orthogonal filtration and Goswami's
quantum isometry group of an admissible spectral triple.
\end{abstract}

\section*{Introduction}
The quantum isometry group of a noncommutative Riemannian compact manifold (an admissible spectral triple)
was defined and constructed by Goswami in~\cite{gos_qiso}.
His breakthrough construction, technically more involved than the previous approaches to quantum symmetry groups
in the case of finite structures~\cite{wang_symgroups, BBjaco}, provides a very natural direct link
between Connes' noncommutative geometry~\cite{ncg}
and the theory of compact quantum groups introduced by Woronowicz in the eighties~\cite{wor_pseudo}.
We refer the reader to the introduction and bibliography of~\cite{QSG_ortho_filtration} for an overview
of the several developments since Goswami's paper.

Motivated by the work of Goswami, Banica and Skalski define and construct in~\cite{QSG_ortho_filtration}
the quantum symmetry group of a $C^*$-algebra endowed with an orthogonal filtration.
Their construction provides a general powerful tool to define and check the existence of quantum symmetry groups
of various mathematical systems and unifies several known quantum symmetry groups constructions.
The work of Banica and Skalski also has the merit to clearly exhibit some of the structures needed to enable one
to prove the existence of a compact quantum symmetry group, see~\cite{QSG_ortho_filtration} for details.
However, although Goswami's work was one of the inspirations for~\cite{QSG_ortho_filtration},
it seems that Goswami's quantum isometry group in~\cite{gos_qiso} cannot, in general,
be seen as a particular case of the quantum symmetry groups defined in~\cite{QSG_ortho_filtration}
(because the subspace spanned by the eigenvalues of Goswami's Laplacian does not seem to form a subalgebra in general).  

It is the purpose of the present paper to propose a construction that simultaneously generalizes the quantum symmetry groups
of Goswami and of Banica-Skalski. We define and construct the quantum symmetry group of a Hilbert module
endowed with an orthogonal filtration. The concept of Hilbert module endowed with an orthogonal filtration
is inspired by Banica-Skalski's notion of $C^*$-algebra equipped with an orthogonal filtration, and is a natural generalization of it.
Also, to an admissible spectral triple in the sense of~\cite{gos_qiso}, one can associate an appropriate Hilbert module endowed
with an orthogonal filtration, and our quantum symmetry group coincides with the quantum isometry group in~\cite{gos_qiso}.

The concept also has the interest to provide an alternative approach to the quantum isometry group of a spectral triple.
The main difference with the approach of~\cite{gos_qiso} is that,
instead of extracting from the spectral triple an analogue of the Laplacian on functions
(the so-called ``noncommutative Laplacian'') and making appropriate assumptions on its spectrum,
we directly use the Dirac operator of the spectral triple, its spectrum and its natural domain.
We then add assumptions to these data to get the desired orthogonally filtered Hilbert module.
In the case of ordinary compact Riemannian manifolds,
our exterior forms based quantum isometry group coincides with the one of Goswami in~\cite{gos_qiso}.

The paper is organized as follows.
In the first part, we briefly recall some basic definitions about compact quantum groups.
Then we introduce in part 2 the concept of Hilbert module endowed with an orthogonal filtration, 
and  define the category of ``quantum transformation groups'' for a Hilbert module
equipped with an orthogonal filtration (our starting point being the notion of action of a compact quantum group on a Hilbert
module given in~\cite{mod_coaction}). Part 3 is devoted to the proof of the existence of a universal object in this category.
In the last part we discuss some examples and compare our construction with the ones of Goswami and Banica-Skalski
mentioned previously.

\vspace{1em}
\textbf{Notations and conventions} $-$
By \emph{algebra} we will always mean \emph{unital algebra}. So that algebra morphisms are assumed to preserve the units.
The symbol  $\odot$ will denote the algebraic tensor product, while $\otimes$ will denote tensor product of maps,
spatial tensor product of $C^*$-algebras, or exterior tensor product of Hilbert modules.


\section{Compact quantum groups}

We recall here some basic definitions on compact quantum groups.
See~\cite{wor_pseudo, wor_CQG, wang_freeprod} and~\cite{pod_coactions} for more details.

\begin{definition}
A \emph{Woronowicz $C^*$-algebra} is a couple $(Q, \Delta)$, where $Q$ is a $C^*$-algebra
and   $\Delta : Q \rightarrow Q \otimes Q$ is a $*$-morphism such that:
\begin{itemize}
	\item $(\Delta \otimes id_Q)\circ \Delta = (id_Q \otimes \Delta) \circ \Delta$,
	\item the spaces $\vect \{\Delta(Q).(Q \otimes 1_Q) \}$
		and $\vect \{ \Delta(Q).(1_Q \otimes Q)\}$ are both dense in $Q \otimes Q$.
\end{itemize}
\end{definition}

\begin{definition}
Let $(Q_0, \Delta_0)$ and $(Q_1, \Delta_1)$ be Woronowicz $C^*$-algebras.
A \emph{morphism of Woronowicz $C^*$-algebras} from $Q_0$ to $Q_1$ is a $*$-morphism: 
\[
	\mu : Q_0 \rightarrow  Q_1 \quad \mbox{such that} \quad (\mu \otimes \mu) \circ \Delta_0 = \Delta_1 \circ \mu .
\] 
\end{definition}

The category of compact quantum groups is then defined to be the opposite category of the category of Woronowicz $C^*$-algebras.

\begin{definitions}
Let $Q =(Q, \Delta)$ be a Woronowicz $C^*$-algebra.
\begin{itemize}
	\item A \emph{Woronowicz $C^*$-ideal} of $Q$ is a $C^*$-ideal $I$ of $Q$
		such that $\Delta(I) \subset \ker(\pi \otimes \pi),$
		where $\pi : Q \rightarrow Q / I$ is the canonical quotient map.
	\item  A \emph{Woronowicz $C^*$-subalgebra} of $Q$ is a $C^*$-subalgebra $Q'$ of $Q$ such that
		$\Delta(Q') \subset Q'\otimes Q'$. 
\end{itemize}
\end{definitions}

\begin{definition}
Let $Q$ be a Woronowicz $C^*$-algebra. A matrix $(v_{ij})_{1 \leqslant i,j \leqslant n} \in \cM_n(Q)$ is called \emph{multiplicative}
if we have $\Delta(v_{ij}) = \sum\limits_{k = 1}^n v_{ik} \otimes v_{kj}$ for all $i,j$.
\end{definition}

The concept of an action of a quantum group on a $C^*$-algebra is formalized as follows.

\begin{definition}
Let $Q$ be a Woronowicz $C^*$-algebra and let $A$ be a $C^*$-algebra.
A \emph{coaction of $Q$ on $A$} is a $*$-morphism $\alpha : A \rightarrow A \otimes Q$ satisfying:
\begin{itemize}
	\item $(\alpha \otimes id_Q) \circ \alpha = (id_A\otimes \Delta) \circ \alpha$,
	\item $\vect \{\alpha(A).(1 \otimes Q) \}$ is dense in $A \otimes Q$.
\end{itemize}
We say that a coaction $\alpha$ of $Q$ on $A$ is \emph{faithful} if there exists no nontrivial  Woronowicz $C^*$-subalgebra
$Q'$ of $Q$ such that $\alpha(A) \subset A \otimes Q'$. Furthermore if $\tau$ is a continuous linear functional on $A$,
we say that $\alpha$ preserves $\tau$ if $(\tau \otimes id_Q ) \circ \alpha = \tau(\cdot)1_Q$.
\end{definition}


\section{Quantum groups actions on Hilbert modules}

We recall now the definition of an action of  a compact quantum group on a Hilbert module (see~\cite{lance} for background material
on Hilbert modules). Then we introduce the notion of orthogonal filtration on a Hilbert module, give some natural examples
of such objects,  and define what we mean by \textit{preserving the filtration} for an action of a compact quantum  group
on a Hilbert module endowed with an orthogonal filtration.

\begin{definition}
Let $A$ be a $C^*$-algebra. A \emph{(right) pre-Hilbert $A$-module} is a vector space $E$, equipped
with a (right) $A$-module structure together with an $A$-valued inner product $\l \cdot | \cdot \r_A$, that is to say:
\begin{itemize}
	 \item $\forall \xi, \eta, \zeta \in E, \forall a,b \in A$, $\l \xi | \eta a +\zeta b\r_A = \l \xi | \eta \r_A a +\l \xi | \zeta \r_A b ,$
	 \item $\forall \xi, \eta \in E, \l \xi | \eta \r_A^* = \l \eta |\xi \r_A ,$
	 \item $\forall \xi \in E$, $\l \xi | \xi \r_A \geqslant 0$ and if $\l \xi  | \xi \r_A = 0$  then $\xi= 0$.
\end{itemize}
We define a norm $\| \cdot \|_A$  on $E$ by setting for $\xi \in E$, $\| \xi \|_A = {\|\l \xi | \xi \r_A\|}^\frac{1}{2}.$
If furthermore $E$ is complete with respect to this norm, we say that $E$ is a \emph{(right) Hilbert $A$-module}.

We say that $E$ is \emph{full} if the space $\l E |E \r_A = \vect \{ \l \xi | \eta \r_A \tq \xi , \eta \in E \}$ is dense in $A$.
\end{definition}

\textit{Left Hilbert $A$-modules} are defined analogously, except that the $A$-valued inner product
$_A \l \cdot | \cdot \r$ has to be linear in the first variable and antilinear in the second one.
In what follows we will mostly consider right Hilbert modules. Of course, the construction can be adapted for left Hilbert modules. 
\vspace{0.1em} 

The notion of coaction on a Hilbert module is due to Baaj and Skandalis~\cite[Definition~2.2]{mod_coaction}.
But working with Woronowicz $C^*$-algebras instead of Hopf $C^*$-algebras simplifies the original definition:

\begin{definition}
Let $A$ be a $C^*$-algebra and let $E$ be a Hilbert $A$-module.
A \emph{coaction} of a Woronowicz $C^*$-algebra $Q$ on $E$ consists of:
\begin{itemize}
	\item a coaction $\alpha:A \rightarrow A \otimes Q$,
	\item a linear map $\beta : E \rightarrow E \otimes Q$ satisfying:
	\begin{enumerate}
	 	\item $\vect \{\beta (E).(1 \otimes Q)\}$ is dense in $E \otimes Q$,
	 	\item $(\beta \otimes id_Q) \circ \beta = (id_E \otimes \Delta) \circ \beta$,
	 	\item $\forall \xi , \eta \in E , \l \beta(\xi) | \beta(\eta) \r_{A\otimes Q} = \alpha(\l \xi | \eta \r_A)$,
		\item $\forall \xi \in E, \forall a \in A, \beta (\xi.a) = \beta(\xi).\alpha(a)$.
	\end{enumerate}
\end{itemize}
We say that the coaction $(\alpha, \beta)$ of $Q$ on $E$ is faithful if there exists no nontrivial  Woronowicz
$C^*$-subalgebra $Q'$ of $Q$ such that $\beta(E) \subset E \otimes Q'$ (note that we do not require $\alpha$ to be faithful).
\end{definition}

\begin{rmq}
If $(\alpha, \beta)$ is a coaction of a Woronowicz $C^*$-algebra $Q$
on a Hilbert $A$-module $E$, then $\beta : E \rightarrow E\otimes Q$ is necessarily continuous. Indeed:

For all $\xi \in E,$ $ \|\beta(\xi)\|^2_{A\otimes Q} = \|\l \beta(\xi) | \beta(\xi) \r_{A \otimes Q} \|
= \| \alpha (\l \xi | \xi \r_A ) \| \leqslant \|  \l \xi | \xi \r_A \| = \|\xi\|^2_A.$
\end{rmq}

\begin{definition}
Let $A$ be a $C^*$-algebra, let $\tau$ be a faithful state on $A$ and let $E$ be a Hilbert $A$-module.
An \textit{orthogonal filtration} $(\tau, (V_i)_{i \in \cI}, J, \xi_0)$ of $E$ consists of:
\begin{itemize}
	\item a family $(V_i)_{i \in \cI}$ of finite-dimensional subspaces of $E$ such that:
		\begin{enumerate}
		\item for all $i,j \in \cI$ with $ i \neq j$, $\forall \xi \in V_i$ and $\forall \eta \in V_j$,
		$\tau(\l \xi | \eta \r_A) = 0$, 
		\item the space $\cE_0 = \sum\limits_{i \in \cI} V_i$ is dense in $(E, \|\cdot \|_A)$,
		\end{enumerate}
	\item an element $\xi_0 \in E$,
	\item a one-to-one antilinear operator $J : \cE_0 \rightarrow \cE_0$.
\end{itemize}
\end{definition}

\begin{exemples}\label{ex_filtration}
\begin{enumerate}[label=(\arabic*)]
	\item Let $M$ be a compact Riemannian manifold. 
		The space of continuous sections of the bundle of exterior forms on $M$, $ \Gamma (  \Lambda^* M)$,
		is a Hilbert $C(M)$-module. We can equip it with an orthogonal filtration by taking
		 $\tau = \displaystyle\int \cdot\, \d vol$ (where $\d vol$ denotes the Riemannian density of $M$),
		$\xi_0 = m \mapsto 1_{\Lambda^*_m M}$,
		$J :  \Gamma (  \Lambda^* M) \rightarrow \Gamma (  \Lambda^* M)$ the canonical involution
		and  $(V_i)_{i \in \NN}$ the family of eigenspaces of the de Rham operator $D = \overline{\d+\d^*}$.

	\item We recall from~\cite{QSG_ortho_filtration} the definition of a $C^*$-algebra equipped with an orthogonal filtration:
		\begin{definition}
		Let $A$ be a $C^*$-algebra, $\tau$ be a faithful state on $A$ and $(V_i)_{i \in \cI}$ be a family of finite-dimensional
		subspaces of $A$ (with the index set $\cI$ containing a distinguished element $0$).
		We say that $(\tau, (V_i)_{i\in \cI})$ is an \emph{orthogonal filtration of $A$} if:
		\begin{enumerate}
			\item $V_0 = \CC.1_A$,
			\item $\forall i,j \in \cI$ such that $i \neq j$, $\forall a \in V_i$ and $\forall b \in V_j$, $\tau (a^*b) = 0$,
			\item the space $\cA_0 = \sum\limits_{i \in \cI} V_i$ is a dense $*$-subalgebra of $A$.
		\end{enumerate}
		Setting $E = A$ (with its canonical Hilbert $A$-module structure), $\xi_0 = 1_A$ and $J = a \mapsto a^*$,
		then $(\tau, (V_i)_{i \in \cI}, J,\xi_0)$ is an orthogonal filtration of $E$.
		\end{definition}

	\item Let $(\cA, \cH ,D)$ be an  admissible spectral triple in the sense of~\cite{gos_qiso}. We set:
		\begin{enumerate}
			\item $E = A = \overline{\,\cA\,}^{\cL(\cH)}$, 
			\item $\tau =\left\{
				\begin{array}{l}
					a \mapsto \dfrac{Tr_\omega (a|D|^{-p})}{Tr_\omega (|D|^{-p})}
					\text{ if }\cH \text{ is infinite dimensional,} \\ 
					\text{the usual trace otherwise,}
				\end{array}\right.$\\
				where $Tr_\omega$ denotes the Dixmier trace and $p$ is the metric dimension of  $(\cA, \cH ,D)$,
			\item the $(V_i)_{i \in \NN}$ are the eigenspaces of the `noncommutative Laplacian', 
			\item $\xi_0$ and $J$ are respectively the unit and the involution of $A$.
		\end{enumerate}
		The couple $(\tau ,(V_i)_{i\in\NN}$) is not in general an orthogonal filtration of $A$ in the sense of~\cite{QSG_ortho_filtration}
		since $\displaystyle\sum_{i \in \NN} V_i$ is not necessarily a $*$-subalgebra of $A$.
		However, $(\tau ,(V_i)_{i\in\NN},J , \xi_0)$ is an orthogonal filtration of $A$, seen as a Hilbert $A$-module.

	\item Let us recall some common conditions on spectral triples.
	\begin{definition}
	Let $(\cA, \cH , D)$ be a spectral triple with finite metric dimension $p$. 
	\begin{itemize}
 		\item We say that $(\cA, \cH, D)$ \emph{satisfies the finiteness and absolute continuity condition} if
			 the space $\displaystyle \cH^\infty = \bigcap_{k \in \NN} \dom(D^k)$ is a finitely generated projective
			left $\cA$-module, and if there exists $q \in \cM_n (\cA)$ with $q = q^2=q^*$ such that:
			\begin{enumerate}
				\item $\cH^\infty \cong \cA^nq$,
				\item the left $\cA$-scalar product $_\cA\l \cdot | \cdot \r$ induced on $\cH^\infty$
					by the previous isomorphism satisfies:
				\[ 
					\frac{Tr_\omega (_\cA\l \xi | \eta \r |D|^{-p})}{Tr_\omega (|D|^{-p})} = (\eta | \xi )_\cH.
				\]
			\end{enumerate}
			(Note that if $(\cA, \cH, D)$ is regular then $\cH^\infty$ is automatically a left $\cA$-module.)

		\item We say that $(\cA, \cH ,D)$ is \textit{real} if it is equipped with an antiunitary operator $\cJ : \cH \rightarrow \cH$
		such that:
			\begin{enumerate}
				\item $\cJ(\dom(D)) \subset \dom(D)$,
				\item $\cJ^2 = \varepsilon$ and  $\cJ D = \varepsilon' D\cJ$, where $\varepsilon, \varepsilon' \in \{-1,1\}$,
				\item $\forall a,b \in \cA$, $[a, \cJ b^*\cJ^*] = 0$.
			\end{enumerate} 
	\end{itemize}
	If $(\cA ,\cH ,D)$ satisfies the finiteness and absolute continuity condition it is natural to consider
	$A =  \overline{\,\cA\,}^{\cL(\cH)}$ and the Hilbert $A$-module $E$ obtained by completing $\cH^\infty$ (for the $A$-norm). 
	The eigenspaces $(V_i)_{i\in \NN}$ of $D$ are two by two orthogonal in $\cH$, thus
	$\forall \xi \in V_i$, $\forall \eta \in V_j$ such that $i \neq j$, we get $\tau (_A \l \xi | \eta \r) = 0$ 
	where $\tau = a \mapsto \dfrac{Tr_\omega (a|D|^{-p})}{Tr_\omega (|D|^{-p})}$.
	If $\tau$ is faithful and $\cE_0$ is dense in $E$, then $E$ can be equipped with an
	orthogonal filtration (with $J : \cE_0 \rightarrow \cE_0$ any one-to-one antilinear map and e.g. $\xi_0 = 0$).

	\noindent If we assume furthermore that $(\cA ,\cH ,D)$ is real, then a natural choice is to set $J = \cJ|_{\cE_0}$.
	\end{definition}

\end{enumerate}
\end{exemples}

\begin{notation}
Let $A$ be a $C^*$-algebra and let $E$ be a Hilbert $A$-module endowed with an orthogonal filtration $(\tau, (V_i)_{i \in \cI}, J, \xi_0)$.
We define on $E$ a scalar product by:
\[
	\forall \xi, \eta \in E, ( \xi | \eta )_\tau = \tau (\l \xi | \eta \r_A).
\]
We denote by $\cH$ the completion of $E$ with respect to this scalar product and by $\|\cdot \|_\tau$
the norm associated with it.
Remark that $\cE_0 = \bigobot{i \in \cI} V_i \subset \cH$ and since 
$\|\xi \|_\tau^2 = \tau (\l \xi | \xi \r_A) \leqslant \| \l \xi | \xi \r_A \| = \|\xi \|_A^2$ for all $\xi \in E$,
we have a continuous injection $E \hookrightarrow \cH$ with dense image.
\end{notation}

We will define now the coactions that preserve the structure of a given Hilbert module equipped with an
orthogonal filtration. This will allow us to describe the category of  its ``quantum transformation groups''.

\begin{definition}
Let $A$ be a $C^*$-algebra and let $E$ be a Hilbert $A$-module endowed with an orthogonal filtration $(\tau, (V_i)_{i \in \cI}, J, \xi_0)$.
A \textit{filtration preserving coaction} of a Woronowicz $C^*$-algebra $Q$ on $E$ is a coaction
$(\alpha,\beta)$ of $Q$ on $E$ satisfying:
\begin{itemize}
	\item $(\tau \otimes id_\cQ) \circ \alpha = \tau(\cdot)1_Q$,
	\item $\forall i \in \cI , \beta (V_i) \subset V_i \odot Q$,
 	\item $(J \otimes *) \circ \beta = \beta \circ J$ on $\cE_0$, where $*$ denotes the involution of $Q$,
	\item $\beta ( \xi_0) = \xi_0 \otimes 1_Q$.
\end{itemize}
In that case, we will also say that $Q$ coacts on $E$ in a \textit{filtration preserving way.}
\end{definition}

\begin{definition}
Let $A$ be a $C^*$-algebra and let $E$ be a Hilbert $A$-module equipped with an orthogonal filtration 
$(\tau, (V_i)_{i \in \cI}, J, \xi_0)$. We will denote by $\cat$ the category of Woronowicz $C^*$-algebras
coacting on $E$ in a filtration preserving way.
If $(\alpha_0, \beta_0)$ and $(\alpha_1,\beta_1)$ are filtration preserving coactions of Woronowicz $C^*$-algebras
$Q_0$ and $Q_1$ on $E$, then a morphism from $Q_0$ to $Q_1$ in this category is a morphism of Woronowicz
$C^*$-algebras $\mu : Q_0 \rightarrow Q_1$ satisfying:
\[
	\alpha_1 = (id_A \otimes \mu) \circ \alpha_0 \quad \text{and} \quad
	\beta_1 = (id_E \otimes \mu) \circ \beta_0.
\]
\end{definition}

\begin{rmq}
If $E$ is full and $\mu: Q_0 \rightarrow Q_1$ is a morphism of Woronowicz $C^*$-algebras satisfying
$\beta_1 = (id_E \otimes \mu) \circ \beta_0$, then $\mu$ automatically satisfies $\alpha_1 = (id_A \otimes \mu) \circ \alpha_0$.

Indeed, for all $\xi, \eta \in E$:
\begin{align*}
\alpha_1(\langle \xi | \eta \rangle_A) 
	& = \langle \beta_1 (\xi) | \beta_1 ( \eta) \rangle_{A \otimes Q_1}
		= \langle (id_E \otimes \mu) \circ \beta_0 (\xi) | (id_E \otimes \mu) \circ \beta_0 (\eta)\rangle_{A \otimes Q_1} \\
	& = (id_A \otimes \mu)(\langle \beta_0(\xi) | \beta_0(\eta) \rangle_{A\otimes Q_0})
		=(id_A \otimes \mu) \circ \alpha_0 (\langle \xi | \eta \rangle_A).
\end{align*}
And since $E$ is full, we get  $\alpha_1 = (id_A \otimes \mu) \circ \alpha_0$.
\end{rmq}

\begin{rmq}
When $E = \Gamma(\Lambda^*M)$ is equipped with the orthogonal filtration 
$(\tau, (V_i)_{i\in\NN} , J, \xi_0)$ described in Example~\ref{ex_filtration}.(1),
the full subcategory of $\cC(E, \tau, (V_i)_{i\in\NN} , J, \xi_0)$ consisting of the commutative Woronowicz $C^*$-algebras coacting
on $\Gamma( \Lambda^* M)$ in a filtration preserving way is antiequivalent to the
category of compact groups acting isometrically on $M$ (see section~\ref{ex_manifold} for more details).
This explains our choice of seeing the opposite category of $\cat$ as the category of quantum transformation groups of $E$.
Moreover since the isometry group of $M$ is a universal object in the category of compact groups acting isometrically on $M$,
we will define the quantum symmetry group of $E$ as a universal object in $\cat$. Proving the existence of
such a universal object is the aim of the next section.
\end{rmq}


\section{Construction of the quantum symmetry group of a Hilbert module equipped with an orthogonal filtration}

The following theorem generalizes the results of Goswami~\cite{gos_qiso} and Banica-Skalski~\cite{QSG_ortho_filtration}.

\begin{theoreme}\label{th_init}
Let $A$ be a $C^*$-algebra and let $E$ be a \textit{full} Hilbert $A$-module endowed with an orthogonal filtration
$(\tau , (V_i)_{i \in \cI}, J, \xi_0)$. The category $\cat$ admits an initial object, which means that there exists a universal
Woronowicz $C^*$-algebra coacting on $E$ in a filtration preserving way. The quantum group corresponding to that
universal object will be called the \textbf{quantum symmetry group} of $(E, \tau ,(V_i)_{i \in \cI}, J , \xi_0)$.
\end{theoreme}

Examples will be discussed in the next section.
This section is devoted to the proof of Theorem~\ref{th_init}.
The proof mostly consists in carefully adapting Goswami's arguments in~\cite[Section 4]{gos_qiso}.
In what follows $E$ denotes a full Hilbert module over a given $C^*$-algebra $A$,
equipped with an orthogonal filtration $(\tau , (V_i)_{i\in \cI}, J , \xi_0).$

\begin{lemme}\label{corep_unitaire}
Let $( \alpha,\beta)$ be a filtration preserving coaction of a Woronowicz $C^*$-algebra $Q$ on $E$.
The $Q$-linear map $\overline{\beta} : \cE_0 \odot Q \rightarrow \cE_0 \odot Q$ given by
$\overline{\beta} (\xi\otimes x) = \beta(\xi)(1\otimes x)$ extends to a unitary of the Hilbert $Q$-module $\cH \otimes Q$.
\end{lemme}

\begin{preuve}
We have for $\xi, \eta \in \cE_0$ and $x,y \in Q$:
\begin{align*}
 {\l\overline{\beta}(\xi \otimes x) | \overline{\beta}(\eta \otimes y)\r}_Q 
	& = x^* {\langle\beta(\xi)|\beta(\eta)\rangle}_Q \; y \\
	& = x^* (\tau \otimes id)({\langle \beta(\xi)|\beta(\eta)\rangle}_{A\otimes Q})\,y
		&&  \text{(by definition of } ( \cdot | \cdot )_\tau = \tau({\langle\cdot | \cdot\rangle}_A)) \\
	& = x^* (\tau \otimes id) \circ \alpha({\langle\xi | \eta \rangle}_A)\, y \\
	& = x^* \; \tau({\langle \xi |\eta\rangle}_A) \, y
		&& \text{(since } (\tau \otimes id) \circ \alpha = \tau(\cdot)1_Q) \\ 
	& = ( \xi | \eta )_\tau \; x^*y = {\langle \xi \otimes x | \eta \otimes y \rangle}_Q .
\end{align*}
In particular $\overline{\beta}$ is isometric and thus extends to a $Q$-linear isometric operator still denoted by
$\overline{\beta} : \cH \otimes Q \rightarrow \cH \otimes Q$. To show that $\overline{\beta}$ is unitary,
it is enough to check that $\overline{\beta}$ has dense image. Since $\vect \{\beta(E).(1 \otimes Q) \}$ is dense in $E \otimes Q$
and $\cE_0$ is dense in $E$, it follows that $\vect \{ \beta(\cE_0).(1 \otimes Q) \}$ is dense in $E \otimes Q$.
Moreover the canonical injection $E\otimes Q \hookrightarrow \cH \otimes Q$ has dense image, 
so that $\vect \{\beta(\cE_0).(1 \otimes Q)\}$ is also dense in $\cH \otimes Q$.
\end{preuve}

\begin{notation}
We define on $\cE_0$ a left scalar product by:
\[
	_\tau( \xi | \eta ) = \tau (\langle J(\xi) | J(\eta) \rangle_A).
\]
For each $i \in \cI$ we set $d_i = \dim (V_i)$ and we fix:
\begin{itemize}
	\item an orthonormal basis $(e_{ij})_{1 \leqslant j \leqslant d_i}$ of $V_i$ for the right scalar product $( \cdot | \cdot )_\tau$,
	\item an orthonormal basis $(f_{ij})_{1 \leqslant j \leqslant d_i}$ of $V_i$ for the left scalar product $_\tau( \cdot | \cdot ) $.
\end{itemize}
We denote by $p^{(i)} \in GL_{d_i} (\CC)$ the change of basis matrix from $(f_{ij})$ to the basis $(e_{ij})$ of $V_i$
and we set $s^{(i)} =  p^{(i)t} \overline{p^{(i)}}$.
\end{notation}

\begin{lemme}\label{unitaire_fidele}
Let $( \alpha,\beta)$ be a filtration preserving coaction of a Woronowicz $C^*$-algebra $Q$ on $E$.
For all  $i \in \cI$, we  denote by $v^{(i)}$ the multiplicative matrix
associated with the  basis $(e_{ij})_{1 \leqslant j \leqslant d_i}$ of the $Q$-comodule~$V_i$
(in other words, $v^{(i)}$ is characterized by:
$\forall j , \beta(e_{ij}) = \sum\limits_{k = 1}^{d_i} e_{ik} \otimes v_{kj}^{(i)}$). 
\begin{itemize}
	\item For all $i \in \cI$, the matrix $v^{(i)} = (v_{kj}^{(i)})_{1\leqslant k,j \leqslant d_i}$ is unitary and
		\[
			v^{(i)t}s^{(i)} \overline{v^{(i)}} {(s^{(i)})}^{-1}
	  		= s^{(i)} \overline{v^{(i)}} {(s^{(i)})}^{-1} v^{(i)t} = I_{d_i}
		\]
	\item The unital $C^*$-subalgebra $Q'$ of $Q$ generated by $\{v_{kj}^{(i)} \tq i \in \cI, j,k \in \{1,\ldots,d_i\}\}$
		is a Woronowicz $C^*$-subalgebra of $Q$ satisfying $\alpha(A) \subset A \otimes Q'$
		and $\beta(E)\subset E \otimes Q'$.\\
		Furthermore $(\alpha ,\beta)$ is a faithful filtration preserving coaction
		of $Q'$ on $E.$
\end{itemize}
\end{lemme}

\begin{preuve}
First let us check that the $v^{(i)}$'s are unitary matrices.

Consider the unitary $\overline{\beta} : \cH \otimes Q \rightarrow \cH \otimes Q$ of the Hilbert $Q$-module $\cH \otimes Q$
constructed in the previous lemma. For all $i,j,m,n$, we have
$\l e_{ij} \otimes 1 | \overline{\beta}(e_{mn} \otimes 1) \r_Q = \delta_{im} v_{jn}^{(i)}
	= \sum\limits_{k=1}^{d_i} \l e_{ik} \otimes v_{jk}^{(i)*} | e_{mn} \otimes 1\r_Q.$
Thus for all~$i,j$, we have $\overline{\beta}^*(e_{ij} \otimes 1)  =
\sum\limits_{k=1}^{d_i} e_{ik} \otimes v_{jk}^{(i)*}$. Then we get:
\[
	e_{ij} \otimes 1 = \overline{\beta} \circ \overline{\beta}^* (e_{ij} \otimes 1)
	= \sum\limits_{k = 1}^{d_i} \beta (e_{ik})(1 \otimes v_{jk}^{(i)*})
	= \sum\limits_{k,l = 1}^{d_i}  e_{il} \otimes v_{lk}^{(i)}v_{jk}^{(i)*},
\]
which shows that for all $l,j \in \{1,\ldots,d_i\}$, $\sum\limits_{k = 1}^{d_i} v_{lk}^{(i)}v_{jk}^{(i)*} = \delta_{lj}$,
i.e. $v^{(i)}v^{(i)*} = I_{d_i}$.

Similarly (using $e_{ij} \otimes 1 = \overline{\beta}^* \circ \overline{\beta} (e_{ij} \otimes 1)$)
we get $v^{(i)*}v^{(i)} = I_{d_i} $. Thus the matrices $v^{(i)}$ are unitary.

$\star$ Let us show now that $v^{(i)t}s^{(i)} \overline{v^{(i)}} {(s^{(i)})}^{-1}
  = s^{(i)} \overline{v^{(i)}} {(s^{(i)})}^{-1} v^{(i)t} = I_{d_i}$.

For $i \in \cI$, since $v^{(i)}$ is a multiplicative and unitary matrix in a Woronowicz $C^*$-algebra,
the matrix $v^{(i)t}$ is invertible in $\cM_{d_i} (Q)$ (\cf \cite{maes_vandaele}), so it is enough
to prove that $v^{(i)t} s^{(i)} \overline{v^{(i)}} {(s^{(i)})}^{-1} = I_{d_i}$.

Using Sweedler's notations, we get for $\xi, \eta \in \cE_0$:
\begin{align*}
_{Q}\langle \beta (\xi) |  \beta (\eta) \rangle
	& = \sum \tau \left(\langle J(\xi\0) | J (\eta\0) \rangle_A \right) \xi\1 \eta\1^*
		= \sum (\tau \otimes id) \left(\langle  J(\xi\0) \otimes \xi\1^* | J (\eta\0) \otimes \eta\1^*\rangle_{A \otimes Q} \right) \\
	& = (\tau \otimes id) \left(\langle (J \otimes *)\circ \beta (\xi) | (J \otimes *)\circ \beta (\eta)\rangle_{A \otimes Q} \right)
		= (\tau \otimes id) \left(\langle \beta \circ J (\xi) |\beta \circ J (\eta) \rangle_{A \otimes Q} \right) \\
	& = (\tau \otimes id) \circ \alpha \left(\langle J(\xi) | J(\eta) \rangle_A\right) 
		= \tau  \left(\langle J(\xi) | J(\eta) \rangle_A\right)1_Q = \,_{\tau}(\xi|\eta)1_Q.
\end{align*}
Moreover, since $\beta(V_i) \subset V_i \odot Q$, there exists $w^{(i)} \in \cM_{d_i} (Q)$
such that $\forall j$, $\beta(f_{ij}) = \sum\limits_{k=1}^{d_i} f_{ik} \otimes w_{kj}^{(i)}$.\\
Then we get $_{Q}\langle \beta (f_{ij}) |  \beta (f_{ik}) \rangle
	= \sum\limits_{l,m = 1}^{d_i} {}_{\tau}(f_{il} | f_{im}) w_{lj}^{(i)}w_{mk}^{(i)*}
	= \sum\limits_{l = 1}^{d_i}  w_{lj}^{(i)}w_{lk}^{(i)*} = \,_{\tau} (f_{ij}|f_{ik})1_Q = \delta_{jk} $.\\ 
This shows: 
\begin{equation}\label{w_unit}
	w^{(i)t} \overline{w^{(i)}} = I_{d_i}  .       
\end{equation}
Furthermore we have for all $j$, $e_{ij} = \sum\limits_{k=1}^{d_i} p_{kj}^{(i)} f_{ik}$, thus
$w^{(i)} =  p^{(i)} v^{(i)}  (p^{(i)})^{-1}  $.  Then replacing $w^{(i)}$  in the equality~(\ref{w_unit})
we get $ {(p^{(i)t})}^{-1}   v^{(i)t} p^{(i)t}\overline{p^{(i)}} \overline{v^{(i)}} \overline{{(p^{(i)})}^{-1}} = I_{d_i}$,
which shows that $ v^{(i)t} s^{(i)}\overline{v^{(i)}} {(s^{(i)})}^{-1} = I_{d_i}$, where $s^{(i)} =  p^{(i)t}\overline{p^{(i)}}.$ 

$\star$ It then follows easily that $Q'$ is a Woronowicz $C^*$-subalgebra of $Q$ satisfying $\alpha(A) \subset A \otimes Q'$
and $\beta(E) \subset E \otimes Q'$. Indeed, since the $v^{(i)}$'s are multiplicative matrices, we have
$\Delta (Q') \subset Q' \otimes Q'$ so that $Q'$ is a Woronowicz $C^*$-subalgebra of $Q$.
Moreover $\beta(\cE_0) \subset \cE_0 \odot Q'$, thus $\beta(E) \subset E \otimes Q'$. 
Then for all $\xi, \eta \in \cE_0$, $\alpha(\langle \xi | \eta \rangle_A) = \langle \beta(\xi) | \beta(\eta) \rangle_{A \otimes Q}
\subset \langle E \otimes Q' | E\otimes Q' \rangle_{A \otimes Q} \subset A \otimes Q'.$
This shows $\alpha(A) \subset A \otimes Q'$ since $\langle\cE_0|\cE_0 \rangle_A$ is dense in $A$.

$\star$ It remains to check that $(\alpha, \beta)$  is a faithful filtration preserving coaction of $Q'$ on $E$.

 We only show that $\vect \{\alpha(A).(1 \otimes Q') \}$ and
$\vect \{\beta(E).(1 \otimes Q') \}$ are respectively dense in $A \otimes Q'$ and $E\otimes Q'$
(the other conditions that must satisfy $(\alpha, \beta)$ to be a filtration preserving coaction of $Q'$ on $E$
directly follow from the fact that it is a filtration preserving coaction of $Q$ on $E$).

We have for all $i,j$:
\begin{align*}
\sum\limits_{k=1}^{d_i} \beta(e_{ik})(1 \otimes v_{jk}^{(i)*})
	& =\sum\limits_{k,l=1}^{d_i} (e_{il} \otimes v_{lk}^{(i)})(1 \otimes v_{jk}^{(i)*}) 
		= \sum\limits_{k,l=1}^{d_i} e_{il} \otimes (v_{lk}^{(i)} v_{jk}^{(i)*}) \\
	& = \sum\limits_{l=1}^{d_i} e_{il} \otimes \delta_{lj}
		\qquad (\text{since } v^{(i)} v^{(i)*} = I_{d_i}) \\
	& = e_{ij} \otimes 1.
\end{align*}
This implies that $\vect \{\beta(E).(1\otimes Q')\}$ is dense in $E\otimes Q'$.
Moreover, we have seen that for all $i\in \cI$, $\overline{v^{(i)}}$ is invertible in $\cM_{d_i}(Q')$
with inverse $x^{(i)} =  {(s^{(i)})}^{-1}v^{(i)t}s^{(i)} .$ 

Let $i,j \in \cI$, $m \in \{1,\ldots,d_i\}$ and $n\in \{1,\ldots, d_j\}$:
\begin{align*}\displaybreak[0]
\sum\limits_{k=1}^{d_i}\sum\limits_{l=1}^{d_j}
		\alpha {(\langle e_{ik} |e_{jl}\rangle}_A).\left(1 \otimes v_{nl}^{(j)*}x_{km}^{(i)}\right)
	& = \sum\limits_{k=1}^{d_i}\sum\limits_{l=1}^{d_j} {\langle\beta (e_{ik}) |\beta(e_{jl})\rangle}_{A \otimes Q}.
		\left(1 \otimes v_{nl}^{(j)*}x_{km}^{(i)}\right) \\
	& = \sum\limits_{k=1}^{d_i} \langle \beta(e_{ik}) |\underbrace{\sum\limits_{l=1}^{d_j}
		\beta(e_{jl}).(1 \otimes v_{nl}^{(j)*})}_{e_{jn} \otimes 1}\rangle_{A \otimes Q}.(1 \otimes x_{km}^{(i)}) \\
	& = \sum\limits_{k,l=1}^{d_i} {\langle e_{il} \otimes v_{lk}^{(i)} |
		e_{jn} \otimes 1\rangle}_{A \otimes Q}.(1 \otimes x_{km}^{(i)}) \\
	& = \sum\limits_{k,l=1}^{d_i} {\langle e_{il} | e_{jn}\rangle}_A \otimes (v_{lk}^{(i)*}x_{km}^{(i)}) \\
	& = \sum\limits_{l=1}^{d_i} {\langle e_{il} | e_{jn}\rangle}_A \otimes \delta_{lm} \qquad (\text{car } \overline{v^{(i)}}x^{(i)}
		= I_{d_i}) \\
	& = {\langle e_{im} | e_{jn}\rangle}_A \otimes 1.
\end{align*}
Thus for all $i,j,m,n$, ${\langle e_{im} | e_{jn}\rangle}_A \otimes 1$ is in $\vect\{\alpha(A).(1 \otimes Q')\}$. 
By density of ${\langle\cE_0 | \cE_0\rangle}_A$ in $A$, this shows that $\vect \{\alpha(A).(1 \otimes Q') \}$ is dense in $A \otimes Q'$.
\end{preuve}

\begin{notation}
For all $i \in \cI$, we consider $\cA_u(s^{(i)})$ the universal Woronowicz $C^*$-algebra 
of Van Daele and  Wang (see \cite{wang_vandaele}) associated with $s^{(i)}$.
That is, $\cA_u(s^{(i)})$ is the universal Woronowicz $C^*$-algebra  generated by a multiplicative and unitary matrix
$u^{(i)} = (u_{kj}^{(i)})_{1 \leqslant k,j \leqslant d_i}$, satisfying the following relations:
\[
	u^{(i)t} s^{(i)}\overline{u^{(i)}} {(s^{(i)})}^{-1}
	=  s^{(i)}\overline{u^{(i)}}{(s^{(i)})}^{-1}  u^{(i)t} = I_{d_i}.
\]
We set $\cU = \Ast{i \in \cI} \cA_u(s^{(i)})$ and $\beta_u : \cE_0 \rightarrow \cE_0 \odot \cU$ the linear map given by:
\[
	\beta_u (e_{ij}) = \sum\limits_{k=1}^{d_i} e_{ik} \otimes u_{kj}^{(i)}.
\]
See \cite{wang_freeprod} for the construction of free product of compact quantum groups.

In the following, if $(Q, \Delta)$ is a Woronowicz $C^*$-algebra and $I$ is a Woronowicz
$C^*$-ideal of $Q$, we will denote by $\pi_I : Q \rightarrow Q/I$  the canonical projection and
by $\Delta_I$ the canonical coproduct of $Q/I$ (\ie $\Delta_I : Q/I \rightarrow Q/I \otimes Q/I$
is the unique $*$-morphism satisfying $\Delta_I \circ \pi_I = (\pi_I \otimes \pi_I) \circ \Delta$).
\end{notation}

\begin{lemme}\label{isomorphe_quotient}
Let $( \alpha, \beta)$ be a faithful filtration preserving coaction  of a Woronowicz  $C^*$-algebra $Q$ on $E$.
There exists a Woronowicz $C^*$-ideal $I \subset \cU$ and a faithful filtration preserving coaction
$(\alpha_I , \beta_I)$ of $\cU/I$ on $E$ such that:
\begin{itemize}
	\item $\cU/I$ and $Q$ are isomorphic in $\cat$,
	\item $\beta_I$ extends $(id \otimes \pi_I) \circ \beta_u.$
\end{itemize}
\end{lemme}

\begin{preuve}
For all $i \in \cI$, we denote by $v^{(i)}$ the multiplicative matrix associated with the  basis
$(e_{ij})_{1 \leqslant j \leqslant d_i}$ of the $Q$-comodule~$V_i$. 
In virtue of Lemma~\ref{unitaire_fidele}, we know that $v^{(i)}$ is unitary and satisfies
$v^{(i)t} s^{(i)}\overline{v^{(i)}} {(s^{(i)})}^{-1}
  = s^{(i)} \overline{v^{(i)}} {(s^{(i)})}^{-1} v^{(i)t} = I_{d_i}$.
So by universal property of $\cU$ there exists a morphism of Woronowicz $C^*$-algebras $\mu : \cU \rightarrow Q$ 
such that for all $i,p,q$, $\mu(u_{pq}^{(i)}) = v_{pq}^{(i)}$.
Then $\im\; \mu$ is a Woronowicz $C^*$-subalgebra of $Q$ and for all $i,j$:
\begin{equation}\label{ext}
	\beta (e_{ij}) = \sum\limits_{k=1}^{d_i} e_{ik}\otimes v_{kj}^{(i)} 
	= (id \otimes \mu) \circ \beta_u (e_{ij}).
\end{equation}
Thus the inclusion $\beta(E) \subset E \otimes (\im \; \mu)$ holds, so that $\mu$ is necessarily onto (since  the coaction is faithful).
We set $ I = \ker\ \mu $, we denote by $\widehat{\mu} : \cU/I \rightarrow Q$ the isomorphism of
Woronowicz $C^*$-algebras such that $\widehat{\mu} \circ \pi_I = \mu$ and we set 
$\alpha_I = (id \otimes \widehat{\mu}^{-1}) \circ \alpha$ and  $\beta_I = (id \otimes \widehat{\mu}^{-1}) \circ \beta$.
It is then easy to see that $(\alpha_I, \beta_I)$ is a filtration preserving coaction, and that:
\[
	\widehat{\mu} : (\cU/I,\Delta_I,\alpha_I, \beta_I) \rightarrow (Q,\Delta,\alpha, \beta)
	\;\; \text{is an isomorphism.}
\]
Thanks to~\eqref{ext}, we see that $\beta_I = (id \otimes \widehat{\mu}^{-1}) \circ\beta$ extends
$(id\otimes \pi_I) \circ \beta_u$.
\end{preuve}

Before proving Theorem~\ref{th_init} we need a last lemma.

\begin{lemme}
Let $A$ and $B$ be $C^*$-algebras and let $\I$ be a nonempty family of $C^*$-ideals of $B$. 
Set $I_0 = \bigcap\limits_{I \in \I} I$, and for $I \in \I$,  set  $p_I : B/I_0 \rightarrow B/I$ the unique $*$-morphism such that:
\[
	\xymatrix{
	& B \ar[dl]_*{\pi_0}  \ar[dr]^*{\pi_I}  & \\
	B/I_0 \ar[rr]^{p_I} & & B/I   }
\] 
where $\pi_0$ and $\pi_I$ denote the canonical projections. Then we have for all $x \in A \otimes (B/I_0) $,
\[
	\| x \| = \sup\limits_{I \in \I} \| id \otimes p_I (x) \|.
\]
\end{lemme}

\begin{preuve}
Let $\rho_A : A \hookrightarrow \cL (H)$ be a faithful representation of $A$ on a Hilbert space $H$,
and for $I \in \I$, let $\rho_I : B/I \hookrightarrow \cL(K_I)$ be a faithful representation of $B/I$ on $K_I$.
The $\rho_I$'s induce a faithful representation of $B/I_0$ on $K = \bigoplus\limits_{I\in\I} K_I$:
\[
	\appl[\rho_0:]{B/I_0}{\cL(K)}{x}{\bigoplus\limits_{I \in \I} \rho_I \circ p_I(x),}
\]
We have canonical injections:
\[
	\bigoplus\limits_{I \in \I} \cL(H \otimes K_I) \hookrightarrow  \cL( \bigoplus\limits_{I \in \I} H \otimes K_I) \cong \cL(H \otimes K),
\]
and since injective morphisms of $C^*$-algebras are isometric, we have for $x \in A \otimes B/I_0$:
\begin{align*}
	\|x\| & = \| \rho_A \otimes \rho_0 (x) \| \in \bigoplus\limits_{I \in \I} \cL(H \otimes K_I) \\
		& =\sup\limits_{I \in \I} \| (\rho_A \otimes (\rho_I \circ p_I)) (x) \|  
	 = \sup\limits_{I \in \I} \| (id \otimes p_I)(x)\|. \tag*{\qedhere}
\end{align*}
\end{preuve}

We are now ready to prove Theorem~\ref{th_init}.

\vspace{\topsep}
\begin{preuve}[\ref{th_init}]
We denote by $\I$ the set of all $C^*$-ideals  $I \subset \cU$ such that: 
\begin{center}
$(id \otimes \pi_I) \circ \beta_u$ extends to a continuous linear map
$\beta_I : E \rightarrow E \otimes \cU/I$ such that there exists a $*$-morphism
$\alpha_I : A \rightarrow A \otimes \cU/I$ preserving $\tau$ and satisfying:
\begin{itemize} 
	\item $\forall \xi, \eta \in E, \langle \beta_I(\xi)|\beta_I(\eta) \rangle_{A\otimes \cU/I} = \alpha_I(\langle\xi | \eta \rangle_A)$,
	\item $\forall \xi \in E, \forall a \in A, \beta_I(\xi.a) = \beta_I(\xi).\alpha_I(a)$,
	\item $(J \otimes *) \circ \beta_I = \beta_I \circ J$ on $\cE_0$,
	\item $\beta_I ( \xi_0) = \xi_0 \otimes 1$.
\end{itemize}
\end{center}
The set $\mathscr I$ is nonempty, since it contains the kernel of the counit $\varepsilon : \cU \rightarrow \CC$
(this can be directly checked, or seen by applying Lemma~\ref{isomorphe_quotient} to the trivial coaction
$A \rightarrow A \otimes \CC$, $E\rightarrow E \otimes \CC$).\\
We denote by $I_0$ the intersection of all elements of $\mathscr I$, by $Q_0 = \cU/I_0$ and by $\pi_0 : \cU \rightarrow Q_0$
the canonical projection (as intersection of $C^*$-ideals, $I_0$ is a $C^*$-ideal, so $\pi_0$ is a $*$-morphism). 
Let us show that $I_0 \in \I$.

$\star$ First let us show that $(id \otimes \pi_0) \circ \beta_u$ extends to a continuous linear map
$\beta_0 : E \rightarrow E \otimes Q_0$. Note that for $x \in E \otimes Q_0$,
${\| x \|}_{A\otimes Q_0} = \sup\limits_{I \in \I} {\|id_E \otimes p_I (x) \|}_{A\otimes \cU/I}$.
Indeed:
\begin{align*}
{\| x \|}_{A\otimes Q_0}^2
	& = \| {\langle x | x \rangle}_{A \otimes Q_0} \|
		= \sup\limits_{I \in \I} \| (id_A \otimes p_I)({\langle x | x \rangle}_{A \otimes Q_0}) \| 
		 \qquad(\text{by the previous lemma) }	 \\		
	& = \sup\limits_{I \in \I} \| {\langle id_E \otimes p_I (x) | id_E \otimes p_I(x) \rangle}_{A \otimes \cU/I} \|
		= \sup\limits_{I \in \I} {\|id_E \otimes p_I (x) \|}_{A\otimes \cU/I}^2.
\end{align*}
Furthermore we have for all $\xi \in \cE_0$ and all $I \in \I$:
\[
	\|(id_E \otimes \pi_I)\circ \beta_u(\xi) \|^2_{A \otimes \cU/I} = \| \langle \beta_I (\xi) | \beta_I (\xi) \rangle_{A \otimes \cU/I} \|
	= \| \alpha_I(\langle \xi |\xi\rangle_A)\| \leqslant \|\langle \xi |\xi \rangle_A \| = \|\xi\|^2_{A}
\]
since $\alpha_I$ is a $*$-morphism. Hence for all $\xi \in \cE_0$:
\begin{align*}
\|(id_E \otimes \pi_0) \circ \beta_u(\xi)\|_{A \otimes Q_0}
	& = \sup\limits_{I\in \I} {\|(id_E \otimes (p_I \circ\pi_0)) \circ \beta_u(\xi) \|}_{A\otimes \cU/I} \\
	& = \sup\limits_{I\in \I} {\|(id_E \otimes \pi_I) \circ \beta_u(\xi) \|}_{A\otimes \cU/I} \leqslant \|\xi\|_A,
\end{align*}
which shows that $(id \otimes \pi_0) \circ \beta_u$ extends to a continuous linear map $\beta_0 : E \rightarrow E \otimes Q_0$.

$\star$ Next let us show that there exists a linear map $\alpha_0 : {\langle \cE_0 | \cE_0 \rangle}_A \rightarrow A \otimes Q_0$ 
such that
\[
	\forall \xi, \eta \in \cE_0, \alpha_0({\langle \xi | \eta \rangle}_A) = (id \otimes \pi_0)
	\left({\langle \beta_u (\xi) | \beta_u (\eta) \rangle}_{A\otimes \cU}\right)\! .
\]
Let $\xi_1, \ldots , \xi_n$ and $\eta_1, \ldots \eta_n$ be elements of $\cE_0$
such that $\sum\limits_{i=1}^n {\langle \xi_i | \eta_i \rangle}_A = 0$.

\noindent Then for all $I \in \I$:
\[
	\sum\limits_{i=1}^n (id \otimes \pi_I)\left({\langle \beta_u (\xi_i) | \beta_u (\eta_i) \rangle}_{A\otimes \cU}\right) 
	= \sum\limits_{i=1}^n {\langle \beta_I (\xi_i) | \beta_I (\eta_i) \rangle}_{A\otimes \cU/I}
	= \sum\limits_{i=1}^n \alpha_I ( {\langle \xi_i | \eta_i \rangle}_A ) = \alpha_I(0) = 0.
\]
Thus we have: 
\[
	\left\| \sum\limits_{i=1}^n (id \otimes \pi_0)\left({\langle \beta_u (\xi_i) | \beta_u (\eta_i) \rangle}_{A\otimes \cU}\right) \right\|
	= \sup\limits_{I \in \I} \left\| \sum\limits_{i=1}^n (id \otimes \pi_I)\left({\langle \beta_u (\xi_i) |
	\beta_u (\eta_i) \rangle}_{A\otimes \cU}\right) \right\| = 0.
\]
This shows that we can define a linear map $\alpha_0 : {\langle \cE_0 | \cE_0 \rangle}_A \rightarrow A \otimes Q_0$ by the formula
\[
	\alpha_0 \left( \sum\limits_{i=1}^n {\langle \xi_i | \eta_i \rangle}_A \right)
	= \sum\limits_{i=1}^n (id \otimes \pi_0)\left({\langle \beta_u (\xi_i) | \beta_u (\eta_i) \rangle}_{A\otimes \cU}\right).
\] 

$\star$ Let us check that $\alpha_0 : {\langle \cE_0 | \cE_0 \rangle}_A \rightarrow A \otimes Q_0$ extends to a
$*$-morphism $\alpha_0 : A \rightarrow A \otimes Q_0$ preserving~$\tau$.
We get for all $\xi, \eta \in \cE_0$ and all $I$ in $\I$:
\begin{align*}
(id \otimes p_I) \circ \alpha_0 ({\langle \xi | \eta \rangle}_A)
	&  = (id \otimes p_I) \circ (id \otimes \pi_0)
		\left({\langle \beta_u (\xi) | \beta_u (\eta) \rangle}_{A\otimes \cU}\right) \\
	& = (id \otimes \pi_I) \left({\langle \beta_u (\xi) | \beta_u (\eta) \rangle}_{A\otimes \cU}\right)
		= \alpha_I ({\langle \xi | \eta \rangle}_A).
\end{align*}
Hence $(id \otimes p_I) \circ \alpha_0$ and $\alpha_I$ coincide on ${\langle \cE_0 | \cE_0 \rangle}_A$.
Consequently we have for $x \in {\langle \cE_0 | \cE_0 \rangle}_A$:
\[
	\|\alpha_0(x)\| = \sup\limits_{I \in \I} \| (id \otimes p_I) \circ \alpha_0 (x) \|
	= \sup\limits_{I \in \I} \| \alpha_I (x) \| \leqslant \|x\|.
\]
\nopagebreak Thus $\alpha_0$ extends continuously to $A$.
Moreover for all $a,b \in A$, we have
\begin{align*}
	& (id \otimes p_I) ( \alpha_0 (ab) - \alpha_0(a)\alpha_0(b)) =
		\alpha_I (ab) - \alpha_I(a)\alpha_I(b) = 0 \\
\mbox{and}\quad
	& (id \otimes p_I) (\alpha_0 (a^*) - \alpha_0(a)^*) = \alpha_I (a^*) - \alpha_I(a)^* = 0.
\end{align*}
Hence $\|\alpha_0 (ab) - \alpha_0(a)\alpha_0(b)\| = \sup\limits_{I\in \I} \|(id \otimes p_I) ( \alpha_0 (ab) - \alpha_0(a)\alpha_0(b))\| = 0$.

Similarly, we get  $\|\alpha_0 (a^*) - \alpha_0(a)^*\| = 0.$  So $\alpha_0$ is indeed a $*$-morphism.
Moreover for all $a \in A $ and all $I \in \I$, we have
\begin{align*}
p_I \circ (\tau \otimes id) (\alpha_0(a))
	& =( \tau \otimes p_I)(\alpha_0 (a)) = (\tau \otimes id) \circ \alpha_I (a)\\
	& = \tau(a)1_{\cU/I}  = p_I (\tau(a)1_{Q_0}) .
\end{align*}
Thus $\alpha_0$ preserves $\tau.$ 

$\star$ We are now ready to check that $I_0 \in \I$. For all $\xi, \eta \in \cE_0$
we have (by construction of $\alpha_0$) that
$\alpha_0(\langle\xi | \eta \rangle_A) =\langle \beta_0(\xi)|\beta_0(\eta) \rangle_{A\otimes Q_0}$, 
and this equality extends by continuity for $\xi ,\eta \in E$.
Since $(id_A \otimes p_I)\circ \alpha_0 = \alpha_I$ and $(id_E \otimes p_I)\circ \beta_0 = \beta_I$
for all $I \in \I$, we get for all $\xi \in E$ and $a\in A$:
\begin{align*}
\|\beta_0(\xi .a)-\beta_0(\xi).\alpha_0(a)\|
	& = \sup\limits_{I\in\I}  \| (id_E \otimes p_I) \circ \beta_0 (\xi .a)-(id_E\otimes p_I)\left(\beta_0(\xi).\alpha_0(a)\right) \| \\
	& = \sup\limits_{I\in\I} \| (id_E \otimes p_I) \circ \beta_0 (\xi .a)-[(id_E\otimes p_I)\circ \beta_0 (\xi)].[(id_A\otimes p_I)\circ \alpha_0 (a)]\| \\
	& = \sup\limits_{I\in\I} \| \beta_I(\xi .a)- \beta_I (\xi).\alpha_I (a) \| = 0.
\end{align*}
Thus $\beta_0(\xi .a) = \beta_0(\xi).\alpha_0(a)$.

Similarly for $\xi \in \cE_0$, $ \| (J \otimes *) \circ \beta_0 (\xi) - \beta_0 \circ J (\xi) \|
	= \sup\limits_{I\in\I} \| (J \otimes *) \circ \beta_I (\xi) - \beta_I \circ J (\xi) \| = 0$
and $\| \beta_0(\xi_0) - \xi_0 \otimes 1_{Q_0}\| = \sup\limits_{I\in\I} \|\beta_I(\xi_0) - \xi_0 \otimes 1_{\cU/I}\|=0$.
Thus we have $(J \otimes *) \circ \beta_0 = \beta_0 \circ J$ on $\cE_0$ and $\beta_0(\xi_0) = \xi_0 \otimes 1$.
We conclude that $I_0 \in \I$.

We set $K = \ker ((\pi_0 \otimes \pi_0) \circ \Delta_\cU)$. In order to show that $I_0$ is a Woronowicz $C^*$-ideal we have
to check that $I_0 \subset K$, and by definition of $I_0$ it is enough to show that $K \in \I$. \\
Denote by $\mu : \cU/K \rightarrow \im ((\pi_0 \otimes \pi_0) \circ \Delta_\cU)$ the $C^*$-isomorphism
that satisfies $(\pi_0 \otimes \pi_0) \circ \Delta_\cU = \mu \circ \pi_K$.\\
Then for all $i,j$:
\begin{align}
(id \otimes \mu) \circ  (id \otimes \pi_K) \circ \beta_u (e_{ij})
	& = (id \otimes \pi_0\otimes \pi_0) \circ (id\otimes \Delta_\cU) \circ \beta_u (e_{ij}) \nonumber\\
	& = \sum\limits_{k,l=1}^{d_i} e_{il} \otimes \pi_0(u_{lk}^{(i)}) \otimes \pi_0(u_{kj}^{(i)})
		= \sum\limits_{k=1}^{d_i} (\beta_0 \otimes id) (e_{ik} \otimes \pi_0(u_{kj}^{(i)}))  \nonumber\\
	& = (\beta_0 \otimes id) \circ \beta_0 (e_{ij}).  \label{eq_2}
\end{align}
Thus we have $(\beta_0 \otimes id) \circ \beta_0 (E) \subset E \otimes \im\; \mu$, and we set
$\beta_K = (id \otimes \mu^{-1}) \circ (\beta_0 \otimes id) \circ \beta_0 : E \rightarrow E \otimes \cU/K.$ 
We get for all $i,j,m,n$:
\begin{align}
(\alpha_0 \otimes id)\circ \alpha_0(\langle e_{ij} |e_{mn}\rangle_A) 
	& = \alpha_0\otimes id\left(\langle \sum_{k=1}^{d_i} e_{ik} \otimes \pi_0(u_{kj}^{(i)})
		| \sum_{l=1}^{d_m} e_{ml} \otimes \pi_0(u_{ln}^{(m)}) \rangle_{A\otimes Q_0}\right) \nonumber\\
	& = \alpha_0 \otimes id \left( \sum_{k,l} \langle e_{ik}|e_{ml}\rangle_A \otimes \pi_0(u_{kj}^{(i)*}u_{ln}^{(m)})\right) \nonumber\\
	& = \sum_{k,l} \langle \beta_0(e_{ik}) | \beta_0 (e_{ml}) \rangle_{A \otimes Q_0} \otimes
		\pi_0(u_{kj}^{(i)*}u_{ln}^{(m)}) \nonumber\\\displaybreak[0]
	& = \sum_{k,l} \langle \sum_{p=1}^{d_i} e_{ip} \otimes \pi_0(u_{pk}^{(i)}) | \sum_{q=1}^{d_m} e_{mq} \otimes
		\pi_0(u_{ql}^{(m)})  \rangle_{A \otimes Q_0} \otimes \pi_0(u_{kj}^{(i)*}u_{ln}^{(m)}) \nonumber\\
	& = \sum_{k,l,p,q} \langle e_{ip} |e_{mq} \rangle_A \otimes \pi_0(u_{pk}^{(i)*}u_{ql}^{(m)}) \otimes
		\pi_0(u_{kj}^{(i)*}u_{ln}^{(m)}) \nonumber\\
	& = \langle (\beta_0 \otimes id) \circ \beta_0(e_{ij}) | 
		(\beta_0 \otimes id) \circ \beta_0(e_{ml}) \rangle_{A\otimes Q_0 \otimes Q_0}\nonumber\\
	& = (id \otimes \mu)(\langle \beta_K (e_{ij}) | \beta_K (e_{ml}) \rangle). \label{eq_3}
\end{align}
Hence for $\xi, \eta \in E$, we have
\[
	(\alpha_0 \otimes id) \circ \alpha_0 ({\langle \xi | \eta \rangle}_A) 
	= {\langle (\beta_0 \otimes id) \circ \beta_0 (\xi) | (\beta_0 \otimes id) \circ \beta_0 (\eta) \rangle}_{A\otimes Q_0 \otimes Q_0}
	\in A \otimes \im\; \mu.
\]
Thus we also have $(\alpha_0 \otimes id) \circ \alpha_0 (A) \subset A \otimes \im\; \mu$, and we define:
\[
	\alpha_K = (id \otimes \mu^{-1}) \circ (\alpha_0 \otimes id) \circ \alpha_0 : A \rightarrow A \otimes \cU/K.
\]

We know from~\eqref{eq_2} that $\beta_K$ extends $(id \otimes \pi_K) \circ \beta_u$
and from~\eqref{eq_3} that for all $\xi, \eta \in E$,
$\alpha_K(\langle \xi |\eta \rangle_A) = \langle \beta_K (\xi)|\beta_K(\eta)\rangle_{A\otimes \cU/K}$.

$\star$ Let us check that $\alpha_K$ preserves $\tau$. We have for all $a \in A$:
\begin{align*}
(\tau \otimes id) \circ \alpha_K (x)
	& = (\tau \otimes id_{\cU/K}) \circ (id_A \otimes \mu^{-1}) \circ (\alpha_0 \otimes id_{Q_0}) \circ \alpha_0 (x) \\
	& = \mu^{-1} \circ (\tau \otimes id_{Q_0} \otimes id_{Q_0}) \circ (\alpha_0 \otimes id_{Q_0}) \circ \alpha_0 (x) \\
	& = \mu^{-1} \circ (\tau(\cdot)1_{Q_0} \otimes id_{Q_0}) \circ \alpha_0 (x)
		= \mu^{-1} \left( 1_{Q_0} \otimes (\tau \otimes id_{Q_0}) \circ \alpha_0 (x) \right) \\
	& =  \mu^{-1} \left( 1_{Q_0} \otimes \tau (x) 1_{Q_0} \right) = \tau(x)1_{\cU/K}.
\end{align*}

$\star$ We have for $\xi  \in E$ and $a \in A$:
\begin{align*}
\beta_K(\xi .a) 
	& = (id \otimes \mu^{-1}) \circ (\beta_0 \otimes id) (\beta_0(\xi .a))
		=  (id \otimes \mu^{-1}) \circ (\beta_0 \otimes id) (\beta_0(\xi).\alpha_0 (a)) \\
	& = (id\otimes \mu^{-1}) \left((\beta_0 \otimes id) \circ \beta_0(\xi).(\alpha_0 \otimes id) \circ\alpha_0 (a) \right)
		= \beta_K(\xi).\alpha_K(a).
\end{align*}

$\star$ Moreover, we have on $\cE_0$:
\begin{align*}
(J \otimes *) \circ \beta_K 
	& = (J \otimes *) \circ (id_E \otimes \mu^{-1}) \circ (\beta_0 \otimes id)  \circ \beta_0
		=  (id_E \otimes \mu^{-1}) \circ (J \otimes * \otimes *) \circ (\beta_0 \otimes id) \circ \beta_0 \\
	& = (id_E\otimes \mu^{-1}) \circ (\beta_0 \otimes id) \circ (J \otimes *) \circ \beta_0
		= (id_E\otimes \mu^{-1}) \circ (\beta_0 \otimes id) \circ  \beta_0 \circ J = \beta_K \circ J.\\
\text{and } \beta_K (\xi_0)
	& = (id_E \otimes \mu^{-1}) \circ (\beta_0 \otimes id)  \circ \beta_0(\xi_0)
		= (id_E \otimes \mu^{-1}) \circ (\beta_0 \otimes id)  (\xi_0 \otimes 1) \\
	& = (id_E \otimes \mu^{-1}) (\xi_0 \otimes 1 \otimes 1) = \xi_0 \otimes 1.
\end{align*}

So $K \in \I$ and $I_0$ is indeed a Woronowicz $C^*$-ideal. We denote by $\Delta_0$ the coproduct on $Q_0$.
In order to show that $Q_0 \in \cat$ it only remains to check that $\alpha_0$ and $\beta_0$ are coassociative and
that $\vect \{\alpha_0(A).(1 \otimes Q_0)\}$ and $\vect \{\beta_0(E).(1 \otimes Q_0) \}$
are respectively dense in $A \otimes Q_0$ and $E \otimes Q_0$.

$\star$ We have seen (\emph{cf.}~\eqref{eq_2}) that for all $i,j$, 
$(\beta_0 \otimes id) \circ \beta_0 (e_{ij}) =
(id \otimes \pi_0\otimes \pi_0)\circ (id \otimes \Delta_\cU) \circ \beta_u (e_{ij})$.
But $(\pi_0 \otimes \pi_0)\circ \Delta_\cU = \Delta_0 \circ \pi_0$. Thus:
\[
	(\beta_0 \otimes id) \circ \beta_0 (e_{ij}) 
	= (id \otimes \Delta_0)\circ (id \otimes \pi_0) \circ \beta_u (e_{ij})
	= (id \otimes \Delta_0)\circ \beta_0 (e_{ij}).
\]
We deduce that $(\beta_0 \otimes id) \circ \beta_0  = (id \otimes \Delta_0)\circ \beta_0 $ on $E$.
Hence for $\xi, \eta \in E$:
\begin{align*}
(\alpha_0 \otimes id_{Q_0}) \circ \alpha_0 ( {\langle \xi | \eta \rangle}_A)
	& = {\langle (\beta_0 \otimes id_{Q_0}) \circ \beta_0 (\xi) 
		| (\beta_0 \otimes id_{Q_0}) \circ \beta_0 (\eta) \rangle}_{A\otimes Q_0 \otimes Q_0} \\
	& = {\langle (id_E \otimes \Delta_0)\circ \beta_0 (\xi)
		| (id_E \otimes \Delta_0)\circ \beta_0 (\eta) \rangle}_{A\otimes Q_0 \otimes Q_0} \\
	& = (id_A \otimes \Delta_0) \left( {\langle \beta_0 (\xi) | \beta_0 (\eta) \rangle}_{A \otimes Q_0} \right) \\
	& = (id_A \otimes \Delta_0) \circ \alpha_0 ({\langle \xi | \eta \rangle}_A),
\end{align*}
which shows (by density of ${\langle E |E \rangle}_A$ in $A$) that $\alpha_0$ is coassociative as well.

$\star$ Finally, to show that $\vect \{\alpha_0(A).(1 \otimes Q_0) \}$ and $\vect \{\beta_0(E).(1 \otimes Q_0) \}$
are respectively dense in $A \otimes Q_0$ and $E \otimes Q_0$, we can proceed in the same way
as in the proof of Lemma~\ref{unitaire_fidele},  by checking that for all $i,j$:
\[
	\sum\limits_{k=1}^{d_i} \beta_0(e_{ik})(1 \otimes \pi_0(u_{jk}^{(i)*}))
	= e_{ij} \otimes 1
\]
and for all $i,j,m,n$:
\[
	\sum\limits_{k=1}^{d_i} \sum\limits_{l=1}^{d_j} \alpha_0({\langle e_{ik} | e_{jl} \rangle}_A)
	.\left(1 \otimes \pi_0 \left( u_{nl}^{(j)*} x_{mk}^{(i)} \right) \right)
	= {\langle e_{im} | e_{jn} \rangle}_A \otimes 1,
\]
where $ x^{(i)} = (s^{(i)})^{-1} u^{(i)t}s^{(i)}$ is the inverse of $\overline{u^{(i)}}$.
Thus $(\alpha_0, \beta_0)$ is a filtration preserving coaction of $Q_0$ on $E$.

It remains to see that $Q_0$ is in fact an initial object in the category $\cat$. 
Let $I \subset \cU$ be a Woronowicz $C^*$-ideal such that there exists a filtration preserving coaction
$(\alpha_I , \beta_I )$ of $\cU/I$ on $E$ such that $\beta_I$ extends $(id \otimes \pi_I) \circ \beta_u$.
We get in particular $I \in \I$, thus $I_0 \subset I$ and $p_I : Q_0 \rightarrow \cU/I$  is then a morphism in $\cat$.
 
Such a morphism is unique. Indeed, if $\eta$ is a morphism from $Q_0$ to $\cU/I$
then $(id \otimes \eta) \circ \beta_0 = \beta_I$, so for all $i,j,k$, $\eta \circ \pi_0 (u^{(k)}_{ij}) = \pi_I (u^{(k)}_{ij})$.
Hence $\eta \circ \pi_0 = \pi_I$, and $\eta = p_I$ follows from uniqueness in the factorization theorem.
Finally, according to Lemmas~\ref{unitaire_fidele} and~\ref{isomorphe_quotient}, we conclude that
$(Q_0, \Delta_0, \alpha_0,\beta_0)$ is an initial object in $\cat$.
\end{preuve}

\begin{rmqs}
As in~\cite{QSG_ortho_filtration}, we can make the following remarks:
\begin{itemize}
	\item If $Q \in \cat$ coacts faithfully on $E$, then the morphism $\mu : Q_0 \rightarrow Q$
		is onto. So that the quantum group associated with $Q$ is a quantum subgroup
		of the one associated with $Q_0$.
	\item If $(W_j)_{j \in \cJ}$ is a subfiltration of $(V_i)_{i \in \cI}$
		(that is $(W_j)_{j \in \cJ}$ is an orthogonal filtration of $E$,
		such that $\forall j \in \cJ$,  there exists $i \in \cI$ such that $W_j \subset V_i$)
		then the quantum symmetry group of $(E, \tau,  (W_j)_{j \in \cJ}, J, \xi_0)$ 
		is a quantum subgroup of the quantum symmetry group of $(E,\tau, (V_i)_{i \in \cI}, J , \xi_0)$.
\end{itemize}
 \end{rmqs}


\section{Examples}

\subsection{Example of a $\boldsymbol{C^*}$-algebra equipped with an orthogonal filtration}

We recall from~\cite{QSG_ortho_filtration} the construction of the quantum symmetry group
of a $C^*$-algebra equipped with an orthogonal filtration.

\begin{definition}
Let  $(A,\tau, (V_i)_{i\in \cI})$ be a  $C^*$-algebra equipped with an orthogonal filtration
 (see Example~\ref{ex_filtration} for the definition).
We say that a Woronowicz $C^*$-algebra $Q$ coacting on $A$ coacts in a filtration preserving way,
if the coaction $\alpha : A \rightarrow A \otimes Q$ of $Q$ on $A$ satisfies
for all $i \in \cI, \alpha(V_i) \subset V_i \odot Q$.
\end{definition}

\begin{theoreme}[\cite{QSG_ortho_filtration}]
Let  $(A,\tau, (V_i)_{i\in \cI})$ be a  $C^*$-algebra equipped with an orthogonal filtration.
The category of Woronowicz $C^*$-algebras coacting on $A$ in  a filtration preserving way
admits an initial object. The quantum group corresponding to that universal object is called the
\textbf{quantum symmetry group} of $(A, \tau, (V_i)_{i\in \cI})$.
\end{theoreme}

Setting $E = A$, $\xi_0 = 1_A$ and $J = a \mapsto a^*$, it is easy to see that the quantum symmetry group of
$(E,\tau, (V_i)_{i \in \cI}, J, \xi_0)$ coincides with the one constructed in the previous theorem
(if $(Q,\Delta,\alpha)$ coacts on $(A,\tau, (V_i)_{i\in \cI})$ in a filtration preserving way
then $(\tau \otimes id) \circ \alpha = \tau(\cdot)1_A$ is automatic since $V_0 = \CC.1_A$).

\vspace{\topsep}
In fact our construction allows to see that the category of
Woronowicz $C^*$-algebras coacting on $(A,\tau, (V_i)_{i\in \cI})$ in  a filtration preserving way
 admits an initial object, even when the assumption
``$\cA_0$ is a $*$-subalgebra of $A$'' is dropped.

In particular, we see that our construction generalizes the one of~\cite{gos_qiso} in
the sense that if $(\cA, \cH ,D)$ is an  admissible spectral triple and if we set:
\begin{itemize}
	\item $E = A = \overline{\,\cA\,}^{\cL(\cH)}$, 
	\item $\tau =\left\{
		\begin{array}{l}
			a \mapsto \dfrac{Tr_\omega (a|D|^{-p})}{Tr_\omega (|D|^{-p})}
			\text{ if }\cH \text{ is infinite dimensional,} \\ 
			\text{the usual trace otherwise,}
		\end{array}\right.$\\
		where $Tr_\omega$ denotes the Dixmier trace and $p$ is the metric dimension of  $(\cA, \cH ,D)$,
	\item the $(V_i)_{i \in \NN}$ are the eigenspaces of the `noncommutative Laplacian',
	\item $\xi_0$ and $J$ are respectively the unit and the involution of $A$,
\end{itemize}\vspace{-\topsep}
then we recover the quantum isometry group of $(\cA, \cH ,D)$ in the sense of~\cite{gos_qiso}.
 
Given a spectral triple, we have seen in Example~\ref{ex_filtration} another way to
attach an Hilbert module equipped with an orthogonal filtration  to it (induced by $D$ instead of the Laplacian).
For an admissible spectral triple $(\cA, \cH ,D)$ satisfying conditions of Example~\ref{ex_filtration}.(4),
the quantum symmetry group of $\cH^\infty$ and the quantum isometry group of $(\cA, \cH ,D)$ in
the sense of Goswami both exist. We do not know if they coincide in that situation.
But in the case of the spectral triple of a Riemannian compact manifold  the question is solved in the next paragraph.

\subsection{Example of the bundle of exterior forms on a Riemannian manifold}\label{ex_manifold}

Let $M$ be a compact Riemannian manifold. Set $A = C(M)$, $\tau = \displaystyle\int \cdot\, \d vol$
where $\d vol$ denotes the Riemannian density of $M$
and set $E = \Gamma (  \Lambda^* M)$ equipped with its canonical Hilbert $C(M)$-module structure.
We denote by $D = \overline{\d + \d^*} : L^2(  \Lambda^* M) \rightarrow L^2(  \Lambda^* M)$ the de Rham operator.
$D$ is self-adjoint and has compact resolvent. So that $\sp(D)$ can be written as:
$\sp(D) = \{\lambda_i \tq i \in \NN\}$, with
$\lim\limits_{i \rightarrow +\infty} |\lambda_i| = +\infty$
and where each $\lambda_i$ is a real eigenvalue of $D$ with finite multiplicity. 
For $ i \in \NN$ we denote by $V_i$ the
subspace associated with $\lambda_i$ and by $d_i$ the dimension of $V_i$.
Note that $V_i \subset \Gamma^\infty ( \Lambda^* M) $, so $V_i \subset E$.

Clearly, the family $(V_i)_{i \in \NN}$ is an orthogonal filtration of $E$, $E$ is full and $\cH = L^2( \Lambda^* M) $.

We denote by $\xi_0 = m \mapsto 1_{\Lambda^*_m M} \in E$ and by $J : E \rightarrow E$ the canonical involution.

\vspace{1.5em}
\textbf{Comparison with the quantum isometry group of $\boldsymbol{M}$ as defined in~\cite{gos_qiso}}

\nopagebreak
Let  $(\alpha,\beta)$ be a filtration preserving coaction
of a Woronowicz $C^*$-algebra $Q$ on $E$. For $\phi$ a state on $Q$, we set $\beta_\phi = (id \otimes \phi) \circ \beta :  L^2(\Lambda^* M) \rightarrow  L^2(\Lambda^* M).$
Since $\beta$ preserves the filtration, $\beta_\phi$ commutes with $D$ on $\cE_0$.
This implies that $\forall k \in \NN$, $\beta_\phi(\dom( D^k)) \subset \dom (D^k)$ and
$\beta_\phi \circ D^k = D^k \circ \beta_\phi$ on $\dom(D^k)$. 
Thus $\beta_\phi (\Gamma^\infty (\Lambda^* M)) \subset \Gamma^\infty (\Lambda^* M)$ and $\beta_\phi$
commutes with $D^2$ on $\Gamma^\infty (\Lambda^* M)$.
Now for $f \in C^\infty (M)$, we have
$\beta_\phi(f) = \beta_\phi(f.\xi_0) = \alpha_\phi(f).\xi_0 = \alpha_\phi(f)$, where $\alpha_\phi = (id \otimes \phi) \circ\alpha$.
Thus $\alpha_\phi (C^\infty (M)) \subset C^\infty (M)$
and $\alpha_\phi$ commutes with $\cL$ (the Laplacian on functions) on $C^\infty (M)$. 
This shows that $ \alpha$ is an isometric coaction of $Q$ on $C(M)$
in the sense of~\cite{gos_qiso}.

Thus we have a forgetful functor $\cF$ from $\cC(\Gamma(\Lambda^*M), \tau, (V_i)_{i \in \NN}, J, \xi_0 )$ to the category of Woronowicz
$C^*$-algebras coacting isometrically on $M$ in the sense of~\cite{gos_qiso} defined by
$\cF(Q,\alpha, \beta) =(Q,\alpha)$.
This functor is in fact an equivalence of categories. To see this, we show that it is fully faithful and essentially surjective.

Let $\alpha$ be an isometric coaction of a Woronowicz $C^*$-algebra $Q$ on $M$, in the sense of Goswami.
It can be seen, along the lines of~\cite{rigidity_II}, that there is a well defined map
$\overline{\alpha} : \Gamma(\Lambda^*M) \rightarrow \Gamma(\Lambda^*M) \otimes Q$ satisfying
$\overline{\alpha}(f_0 \d f_1 \wedge \ldots \wedge \d f_k ) = \alpha(f_0) (\d \otimes id)\circ\alpha (f_1) \ldots  (\d \otimes id)\circ\alpha (f_k)$
for all $f_0 \in C(M)$, $f_1, \ldots , f_k \in C^\infty (M)$, and that $(\alpha, \overline{\alpha})$ is a coaction
of $Q$ on $\Gamma(\Lambda^*M)$. In order to show that for all $k \in \ZZ$, $\overline{\alpha} (V_k) \subset V_k \odot Q$, we check
that $(D \otimes id) \circ \overline{\alpha} = \overline{\alpha} \circ D$ on $\Gamma^\infty(\Lambda^* M)$.
The equality $(\d \otimes id) \circ \overline{\alpha} = \overline{\alpha} \circ \d$ holds by definition of $\overline{\alpha}.$
Now for all $\omega , \eta \in \Gamma^\infty(\Lambda^*M)$, all $x \in Q$ we have
\begin{align*}
\l (\d^* \otimes id) \circ \overline{\alpha}(\omega) | \overline{\alpha}(\eta).(1\otimes x) \r_Q
	& = \l \overline{\alpha}(\omega) | (\d\otimes id) \circ \overline{\alpha}(\eta) \r_Q .x 
		=  \l \overline{\alpha}(\omega) |  \overline{\alpha}(\d \eta) \r_Q .x \\
	& = (\tau \otimes id)( \l \overline{\alpha}(\omega) |  \overline{\alpha}(\d \eta) \r_{A \otimes Q}).x
		= (\omega | \d \eta)_\cH.x = (\d^* \omega | \eta)_\cH.x \\
	& =  \l \overline{\alpha}(\d^*\omega) | \overline{\alpha}(\eta) \r_Q .x 
		= \l \overline{\alpha}(\d^*\omega) | \overline{\alpha}(\eta).(1 \otimes x) \r_Q,
\end{align*}
and since $\vect \{\overline{\alpha}(\Gamma^\infty (\Lambda^* M)).(1 \otimes Q) \}$ is dense in
$\Gamma(\Lambda^*M) \otimes Q$, we have $(\d^* \otimes id) \circ \overline{\alpha} = \overline{\alpha} \circ \d^*$
on $\Gamma^\infty(\Lambda^* M)$. This allows to see that $(\alpha, \overline{\alpha})$ is  a
filtration preserving coaction of $Q$ on $\Gamma(\Lambda^*M)$ and $\cF(Q, \alpha, \overline{\alpha}) = (Q,\alpha)$,
so $\cF$ is essentially surjective.

Let $(\alpha, \beta)$ be a filtration preserving coaction of a Woronowicz $C^*$-algebra $Q$ on $E$.
Let us show by induction on  $k \in \NN$ that $\beta = \overline{\alpha}$ on $\Gamma(\Lambda^k M)$.
The fact that $\beta = \overline{\alpha} = \alpha$ on $C(M)$ is clear.
Let $k\in \NN$ such that for all $l\leqslant k$, $\beta = \overline{\alpha}$ on $\Gamma(\Lambda^l M)$.
Let $f_0$ be in $C(M)$ and $f_1, \ldots , f_{k+1}$ be in $C^\infty (M)$. Then:
\begin{align*}
 \beta ( \d f_1 \wedge  \ldots \wedge  \d f_{k+1})  
       	& = \beta (D(f_1 \d f_2 \wedge \ldots \wedge \d f_{k+1})) -\beta ( \d^*(f_1 \d f_2 \wedge \ldots \wedge \d f_{k+1})) \\
	& = (D \otimes id) \circ \beta (f_1 \d f_2 \wedge \ldots \wedge \d f_{k+1})
		- \overline{\alpha}( \d^*(f_1 \d f_2 \wedge \ldots \wedge \d f_{k+1})) \\
	&  = (D \otimes id)\circ \overline{\alpha} (f_1 \d f_2 \wedge \ldots \wedge \d f_{k+1}) 
		- (\d^* \otimes id) \circ \overline{\alpha}(f_1 \d f_2 \wedge \ldots \wedge \d f_{k+1}) \\
	& = (\d \otimes id) \circ \overline{\alpha}(f_1 \d f_2 \wedge \ldots \wedge \d f_{k+1}) 
		 = \overline{\alpha}(\d(f_1 \d f_2 \wedge \ldots \wedge \d f_{k+1})) \\
	& = \overline{\alpha} ( \d f_1 \wedge  \ldots \wedge  \d f_{k+1}). 
\end{align*}
Thus we have 
$ \beta (f_0 \d f_1 \wedge  \ldots \wedge  \d f_{k+1})  = \alpha(f_0) \overline{\alpha} ( \d f_1 \wedge  \ldots \wedge  \d f_{k+1})  
	=  \overline{\alpha}  (f_0 \d f_1 \wedge  \ldots \wedge  \d f_{k+1})$, which ends the induction.
Now the fact that $\beta$ necessarily coincides with $\overline{\alpha}$ allows to see that a morphism of Woronowicz
$C^*$-algebras $\mu : (Q_0, \alpha_0 , \beta_0) \rightarrow (Q_1, \alpha_1, \beta_1)$, between
Woronowicz $C^*$-algebras coacting in a filtration preserving way on $\Gamma(\Lambda^* M)$, 
that satisfies $(id \otimes \mu) \circ \alpha_0 = \alpha_1 \circ \mu$, automatically
satisfies  $(id \otimes \mu) \circ \beta_0 = \beta_1 \circ \mu$. This shows that $\cF$ is fully faithful,
so $\cF$ is indeed an equivalence of categories. In particular, $\cF$ preserves the initial object, so our
quantum isometry group of $M$ coincides with the one of Goswami.

Note  that in case $(\alpha, \beta)$ is a filtration preserving coaction of $C(G)$ on $\Gamma(\Lambda^* M)$, for
$G$ a given compact group, then there exists an isometric action $\gamma: M \times G \rightarrow M$ of
$G$ on $M$ such that
\[
	\appl[\alpha=]{C(M)}{C(M\times G) \cong C(M)\otimes C(G)}{f}{f \circ \gamma}
\]
and 
\[
	 \appl[\beta =]{\Gamma(\Lambda^*M)}{C(G,\Gamma(\Lambda^*M)) \cong \Gamma(\Lambda^*M) \otimes C(G)}
	{\omega}{(g\mapsto \gamma_g^* (\omega)),}
\]
where $\gamma_g = m \mapsto \gamma(m,g)$ and $\gamma_g^* : \Gamma(\Lambda^* M) \rightarrow \Gamma(\Lambda^*M)$
denotes the pullback by $\gamma_g$.

\subsection{Basic example: free orthogonal quantum groups}

Let $n$ be in $\NN$. We set $A = \CC$, $E = \CC^n$ equipped with its canonical Hilbert space structure,
$\xi_0 = 0$ and $V_0 = \CC^n$.
Let $J: \CC^n \rightarrow \CC^n$ be any invertible antilinear map. We denote by $P$ the matrix
of $J$ in the canonical basis and by $\cA_o (P)$ the universal Woronowicz $C^*$-algebra generated by a multiplicative and unitary matrix
$u = (u_{ij})_{1 \leqslant i,j \leqslant n}$, satisfying the relation $u = P \overline{u}P^{-1}$ (the quantum group 
associated with $\cA_o(P)$ is a so-called \textit{free orthogonal quantum group}, see~\cite{qortho}).
We denote by $\alpha_P: \CC \rightarrow \CC \otimes \cA_o(P)$ the trivial coaction
and by $\beta_P : \CC^n \rightarrow \CC^n \otimes \cA_o(P)$ the linear map given by
$\beta(e_i) = \sum\limits_{k=1}^n e_k \otimes u_{ki}$ where $(e_k)_{1 \leqslant k \leqslant n}$ is the canonical
basis of $\CC^n$.
We can easily check that $(id_\CC ,  (V_0) , J, \xi_0)$ is an orthogonal filtration of $E$ and
that $(\alpha_P , \beta_P)$ is a coaction of $\cA_o(P)$ on $E$.
To see that $(\alpha_P , \beta_P)$ is a filtration preserving coaction, the only nontrivial point 
is to check that $(J \otimes *)\circ \beta_P = \beta_P \circ J.$
We have  for all $i$ in $\{1, \ldots , n\}$:
\begin{align*}
(J \otimes *) \circ \beta_P (e_i) & = \sum_{k=1}^n J(e_k) \otimes u_{ki}^*
	=\sum_{k, l=1}^n P_{lk}e_l \otimes u_{ki}^* = \sum_{l=1}^n e_l \otimes \left( \sum_{k=1}^n P_{lk}u_{ki}^*\right) \\
	& =  \sum_{l=1}^n e_l \otimes \left( \sum_{k=1}^n u_{lk} P_{ki}\right) \qquad (\text{since } P\overline{u} = uP) \\
	& =  \sum_{k=1}^n P_{ki} \beta_P (e_k) = \beta_P \circ J (e_i).
\end{align*}
So $(\alpha_P, \beta_P)$ is a filtration preserving coaction of $\cA_o(P)$ on $\CC^n$. Now we show that
it is a universal object in the category $\cC (E, id_\CC, (V_0), J, \xi_0)$. Let $(\alpha, \beta)$ be a
filtration preserving coaction of a Woronowicz $C^*$-algebra $Q$ on $E$ and let $v = (v_{ij})_{1 \leqslant i,j \leqslant n} \in \cM_n (Q)$
be characterized by $\beta(e_i) = \sum\limits_{k=1}^n e_k \otimes v_{ki}$.
By Lemma~\ref{unitaire_fidele} we already know that $v$ is unitary. Furthermore, by a similar computation to the previous one, 
we see that $(J \otimes *) \circ \beta = \beta \circ J$ leads to the equality $P \overline{v} = v P$.
Thus by universal property of $\cA_o(P)$ we get the existence of a morphism $\mu : \cA_o(P) \rightarrow Q$ 
such that for all $i,j \in \{1,\ldots , n\},$ $ \mu(u_{ij}) = v_{ij}$, which is clearly a morphism in the category
$\cC (E, id_\CC, (V_0), J, \xi_0)$. Consequently the quantum symmetry group of $E$ is the free orthogonal quantum group
associated with $P$.

\subsection{Example built on segments}

Huang~\cite{huang} has constructed examples of faithful actions of non-classical quantum groups on connected metric spaces.
We now examine some of his examples and put them into our framework.

We first describe in this paragraph the Hilbert module endowed with an orthogonal filtration we associate to $d$ disjoint copies of
$[0,1]$. Then we compute its quantum symmetry group, which appears to be the hyperoctahedral quantum group. Finally, we show the existence of a universal object in a certain subcategory
of the one of Woronowicz $C^*$-algebras coacting on this Hilbert module in a filtration preserving way. The quantum group associated
with that universal object might be seen as a quantum isometry group of a certain quotient of $[0,1] \times  \{1,\ldots,  d\}$.

\vspace{1em}
\textbf{Hilbert module associated with  $\boldsymbol{[0,1] \times  \{1,\ldots,  d\}}$}

\noindent We set $I = [0,1]$ and we denote by $\delta_+ : L^2(I) \rightarrow L^2(I)$ the operator $\frac{d}{dx}$ with domain:
\[
	\dom(\delta_+) = \{ f \in H^1(I) \tq f(0)=f(1)=0 \}.
\]
Its adjoint operator is $\delta_- = -\frac{d}{dx}$ with domain $H^1(I)$. 

We define $D_0 : L^2(\Lambda^*(I)) \rightarrow  L^2(\Lambda^*(I)) \cong   L^2(I) \oplus   L^2(I)$ by:
\[
	D_0 = 
	\begin{pmatrix}
	0 & \delta_- \\
	\delta_+ & 0
	\end{pmatrix}.
\]
$D_0$ is a self-adjoint operator with compact resolvent.
It can be checked that the eigenvectors of $D_0$ are the $(\sin(\pi k \dot) , \cos(\pi k \dot))$ with $k \in \ZZ$.
We set $A =  C(I)^d \cong C(I \times  \{1,\ldots, d\})$, $E =(C_0(I) \oplus C(I))^d$. For $i \in \{1,\ldots, d\}$ and $n \in \ZZ$, we denote
by $e_{ni} \in E$ the vector whose components are zero, except the $i$-th one whose value is $(\sin(\pi n \dot) , \cos(\pi n \dot))$.
We set $V_n = \vect\{e_{ni} \tq i \in  \{1,\ldots, d\} \}$, $\tau = \sum\limits_{i=1}^d \int\! \cdot \, \d x_i$,
$\xi_0 = (\underbrace{(0,1),\ldots,(0,1)}_{d \text{ times}})$ and we denote by $J : E \rightarrow E$ the complex conjugation operator.
Then $(E, \tau, (V_n)_{n \in \ZZ} , J , \xi_0)$ is a Hilbert $A$-module equipped with an orthogonal filtration.
To see that this Hilbert module is a good description of $[0,1] \times  \{1,\ldots, d\}$, just remark
that it is obtained  from the spectral triple $(A, H, D)$, where $H =L^2(\Lambda^*(I))^d$ 
 and $D = diag(\underbrace{D_0,\ldots,D_0}_{d \text{ times}})$ (\cf Example~\ref{ex_filtration}.(4) to see how we associate
a Hilbert module endowed with an orthogonal filtration to a suitable spectral triple).

\vspace{1em}
\textbf{Computation of its quantum symmetry group}\nopagebreak

Let $(\alpha, \beta)$ be a filtration preserving coaction of a Woronowicz $C^*$-algebra $Q$ on $E$.
For $n \in \ZZ$, we denote by $v^{(n)} \in \cM_d(Q)$ the unitary matrix characterized by:
\[
	\forall i \in \{1,\ldots, d\},\; \beta(e_{ni}) = \sum\limits_{j = 1}^d e_{nj} \otimes v_{ji}^{(n)}.
\] 
For $i \in \{1,\ldots,d\}$, let $e_i \in A$  denote the vector whose components are zero
except the $i$-th one which equals $1$, and let  $v$ and $w$  respectively denote $v^{(0)}$ and $v^{(1)}$.

We have $\beta(e_{ni}) = \beta(\overline{e_{ni}}) = (J \otimes *) \circ \beta (e_{ni})
	= \sum\limits_{k=1}^d e_{nk} \otimes v_{ki}^{(n)*}$, thus $v_{ki}^{(n)} = v_{ki}^{(n)*}.$
And
\begin{align*}
\alpha(e_i) & = \alpha(\langle e_{0i} | \xi_0 \rangle_A ) = \langle \beta (e_{0i}) | \xi_0 \otimes 1 \rangle_{A \otimes Q}
		= \sum_{k=1}^d \langle e_{0k} | \xi_0 \rangle_A \otimes v_{ki}^* =  \sum_{k=1}^d e_k \otimes v_{ki} \\
	& = \alpha(e_i^2) = \alpha (e_i)^2 =  \sum_{k=1}^d e_k \otimes v_{ki}^2.
\end{align*}
Thus $v_{ki}^2 = v_{ki}$ for all $i,k$. We get
$ \sum\limits_{k=1}^d v_{ik} =  \sum\limits_{k=1}^d v_{ik} v_{ik}^* = 1$ and
$ \sum\limits_{k=1}^d v_{ki} =  \sum\limits_{k=1}^d v_{ki}^* v_{ki} = 1$ since $v$ is unitary. 
We have furthermore
\[\langle e_{ni} | e_{mj} \rangle_A = \delta_{ij} (\sin(\pi n \dot) \sin(\pi m \dot) + \cos (\pi n \dot) \cos(\pi m \dot))e_i
	= \delta_{ij} \cos(\pi (n-m) \dot)e_i .\]
Thus we get
\begin{align*}
\alpha ( \langle e_{ni} | e_{mj} \rangle_A ) & = \delta_{ij} \alpha ( \cos(\pi(n-m)\dot)e_i) 
		= \langle \beta(e_{ni}) | \beta (e_{mj}) \rangle_{A\otimes Q}
		= \sum\limits_{k,l = 1}^d \langle e_{nk} | e_{ml} \rangle_A \otimes v_{ki}^{(n)}v_{lj}^{(m)}  \\
 	& = \sum\limits_{k = 1}^d \cos(\pi (n-m)\dot) e_k \otimes  v_{ki}^{(n)}v_{kj}^{(m)}
\end{align*}
We deduce that:
\begin{itemize}
	\item for all $i \neq j$, $ \alpha(\langle e_{ni} | e_{mj} \rangle_A ) = 0
		= \sum\limits_{k=1}^d  \cos(\pi (n-m)\dot) e_k \otimes  v_{ki}^{(n)}v_{kj}^{(m)}$, so that for all $n,m,k$,  						$v_{ki}^{(n)}v_{kj}^{(m)} = 0$,
		
	\item $\alpha(\langle e_{ni} | e_{mi} \rangle_A) = \alpha(\cos(\pi(n-m)\dot)e_i)
		= \sum\limits_{k=1}^d \cos(\pi (n-m)\dot) e_k \otimes  v_{ki}^{(n)}v_{ki}^{(m)}.$ \\
		Thus for all $i,k$, and all $n,m,n',m'$ such that $|n-m| = |n'-m'|$,
		$v_{ki}^{(n)}v_{ki}^{(m)} = v_{ki}^{(n')}v_{ki}^{(m')}$.
\end{itemize}
Therefore we have for all $i,j,n$, $ v_{ij}^{(n)} = v_{ij}^{(n)} \left( \sum\limits_{k=1}^d v_{ik}\right) = v_{ij}^{(n)}v_{ij}$
and similarly $ v_{ij}^{(n)} =v_{ij} v_{ij}^{(n)}$.
Consequently, $v_{ij}^{(n)} = v_{ij}^{(n)}v_{ij} = v_{ij}v_{ij}^{(-n)} = v_{ij}^{(-n)}$
(since $|n-0| = |0-n|$), thus $v_{ij}^{(n+1)} = v_{ij}^{(n)}v_{ij}^{(-1)} = v_{ij}^{(n)}w_{ij}$.
Then the fact that $w_{ij}^2 = v_{ij}^{(1)} v_{ij}^{(1)} = v_{ij}$ allows to see by immediate induction that
$v_{ij}^{(n)} =  w_{ij}^{r(n)}$ for all $i,j,n$, where $r(n) = 2$ if $n$ is even and $r(n) = 1$ otherwise.

We obtain finally:
\begin{itemize}
	\item  for $j \neq k$, $w_{ij} w_{ik}  = 0$ and $w_{ji}w_{ki} = S(w_{ij})S(w_{ik}) = S(w_{ik}w_{ij}) = 0$
		(here $S$ denotes the antipode of $Q$).
	\item $\sum\limits_{l = 1}^n w_{il}^2 = \sum\limits_{l = 1}^n v_{il} = 1$
		and similarly  $\sum\limits_{l = 1}^n w_{li}^2 = 1$.
\end{itemize}

This leads to the existence of a unique morphism of Woronowicz $C^*$-algebras $\mu : \cA_h (d) \rightarrow Q$ such
that for all $i,j$, $\mu(u_{ij}) = w_{ij}$ (where $\cA_h(d)$ is the hyperoctahedral quantum group and
the $u_{ij}$'s are the canonical generators of $\cA_h(d)$ - see~\cite{ hyperoctahedral, free_wreath}).

In order to conclude that $\cA_h(d)$ is a universal object in the category of Woronowicz $C^*$-algebras coacting
in a filtration preserving way on $E$, it only remains to check that there exists a  filtration preserving coaction $(\alpha,\beta)$
of $\cA_h(d)$ on $E$ such that for all $n\in \ZZ$ and $i \in  \{1,\ldots, d\}$, 
$\beta(e_{ni}) =   \sum\limits_{k=1}^{d} e_{nk} \otimes u_{ki}^{r(n)}$
(so that the morphism $\mu : \cA_h(d) \rightarrow Q$ constructed previously automatically intertwines the coactions).

For $f \in C([0,1])$, we set $p(f) = x \mapsto \frac{1}{2}(f(x) + f(1-x))$ and $q(f) =x \mapsto \frac{1}{2}(f(x)-f(1-x)).$
We define $\alpha : C([0,1])^d \rightarrow C([0,1])^d \otimes \cA_h(d)$ by:
\[
	\alpha\left(\sum\limits_{i=1}^d f_i.e_i\right) = \sum\limits_{i,k = 1}^d p(f_i).e_k \otimes u_{ki}^2 + q(f_i).e_k \otimes u_{ki}
\]
for all $f_1, \ldots, f_d \in C([0,1]).$
Let us check that $\alpha$ is a $*$-morphism. The fact that $\alpha$ preserves the involution is straightforward,
and we have for all $f_1, \ldots, f_d, g_1, \ldots, g_d \in C([0,1])$
\begin{align*}
\alpha\left(\sum\limits_{i=1}^d  f_i.e_i\right) & \alpha\left(\sum\limits_{j=1}^d g_j.e_j\right) \\\displaybreak[0]
		 & = \left(\sum\limits_{i,k = 1}^d p(f_i).e_k \otimes u_{ki}^2 +
			q(f_i).e_k \otimes u_{ki} \right)\left(\sum\limits_{j,l = 1}^d
			p(g_j).e_l \otimes u_{lj}^2 + q(g_j).e_l \otimes u_{lj} \right) \\
		& = \sum\limits_{i,j,k = 1}^d ( p(f_i).e_k \otimes u_{ki}^2 +
			q(f_i).e_k \otimes u_{ki})( p(g_j).e_k \otimes u_{kj}^2
			+ q(g_j).e_k \otimes u_{kj}) \\
		& = \sum\limits_{i,k = 1}^d (p(f_i)p(g_i) + q(f_i)q(g_i)).e_k \otimes u_{ki}^2
			+ (p(f_i)q(g_i) + q(f_i)p(g_j)).e_k \otimes u_{ki}^2
\end{align*}
(the last equality holds since $u_{ki}u_{kj} = 0$ for $i \neq j$ and $u_{ij}^3 = u_{ij}$).

But we have $p(f)p(g) + q(f)q(g) = p(fg)$ and $p(f)q(g) + q(f)p(g) = q(fg)$ for all $f,g \in C([0,1])$.
We thus get
\[
	\alpha\left(\sum\limits_{i=1}^d f_i.e_i\right)\alpha\left(\sum\limits_{j=1}^d g_j.e_j\right)
		= \alpha \left(\sum\limits_{i=1}^d f_ig_i.e_i\right)
		=\alpha\left(\left(\sum\limits_{i=1}^d f_i.e_i\right)\left(\sum\limits_{j=1}^d g_j.e_j\right)\right).
\]
Now we have for $f \in C([0,1])$ and $i \in \{1, \ldots , d\}$
\begin{align*}
(id \otimes \Delta) \circ \alpha (f.e_i) & = \sum\limits_{k = 1}^d
		p(f).e_k \otimes \Delta(u_{ki})^2 + q(f).e_k \otimes \Delta(u_{ki}) \\
	& = \sum\limits_{k,l = 1}^d p(f).e_k \otimes u_{kl}^2 \otimes u_{li}^2
		+ q(f).e_k \otimes u_{kl} \otimes u_{li},
\end{align*}
and since $q\circ p = p \circ q = 0$, $p^2 = p$ and $q^2 = q$, we also have
\begin{align*}
(\alpha \otimes id) \circ \alpha (f.e_i) & = \sum\limits_{l = 1}^d
		\alpha(p(f).e_l) \otimes u_{li}^2 + \alpha(q(f).e_l) \otimes u_{li} \\
	& = \sum\limits_{k,l = 1}^d p(f).e_k \otimes u_{kl}^2 \otimes u_{li}^2
		+ q(f).e_k \otimes u_{kl} \otimes u_{li}.
\end{align*}
Thus $\alpha$ is coassociative.

We have
$p(\cos(\pi n \dot))  =\dfrac{1+(-1)^n}{2}\cos(\pi n \dot)$ and
$q(\cos(\pi n \dot)) = \dfrac{1+(-1)^{n+1}}{2}\cos(\pi n \dot)$,  so that 
$\displaystyle\alpha (\cos(\pi n \dot)e_i) = \sum\limits_{k=1}^d \cos(\pi n \dot) .e_k \otimes u_{ki}^{r(n)}$
(where $r(n)$ is still equal to $2$ when $n$ is even and to $1$ otherwise).
Therefore, we have for all $n \in \ZZ$ and all $i \in \{1,\ldots , d\}$,
\begin{align*}
\sum\limits_{k=1}^d \alpha (\cos (\pi n \dot).e_k)(1\otimes u_{ik}^{r(n)})
		 & = \sum\limits_{k,l= 1 }^d \cos (\pi n \dot).e_l \otimes u_{lk}^{r(n)}u_{ik}^{r(n)}
		=  \sum\limits_{k= 1 }^d \cos (\pi n \dot).e_i \otimes u_{ik}^{2r(n)} \\
	& = \cos (\pi n \dot).e_i \otimes \left(\sum\limits_{k= 1 }^d  u_{ik}^{2} \right)
		= \cos (\pi n \dot).e_i \otimes 1,
\end{align*}
which allows us to see that $\vect\{\alpha(C(I)^d).(1 \otimes \cA_h(d))\}$ is dense in $C(I)^d \otimes \cA_h(d)$,
so that $\alpha$ is a coaction of $\cA_h(d)$ on $C(I)^d$.
Then we define $\beta : E \rightarrow E \otimes \cA_h(d)$ by:
\[
	\beta((f,g).e_i) =   \sum\limits_{k = 1}^{d} (q(f),p(g)).e_k \otimes u_{ki}^2 + (p(f),q(g)).e_k \otimes u_{ki}
\]
for all $f \in C_0([0,1])$, $g \in C([0,1])$ and $i \in \{1,\ldots, d\}$.

We have for $f, f' \in C_0([0,1]), \: g,g' \in C([0,1])$ and $i,j \in  \{1,\ldots, d\}$
\begin{align*}
\langle \beta  ((f,g).e_i   |  \beta((f',g').e_j ) \rangle_{A \otimes \cA_h(d)} \displaybreak[0]
	& = \delta_{ij}   \sum\limits_{k=1}^{d} (q(\overline{f})q(f') + p(\overline{g})p(g') +p(\overline{f})p(f')
		+ q(\overline{g})q(g')).e_k \otimes u_{ki}^2  \\
	&  \qquad \quad \!+(q(\overline{f})p(f') + p(\overline{g})q(g') +p(\overline{f})q(f') + q(\overline{g})p(g')).e_k
		\otimes u_{ki} \\
	& = \delta_{ij}   \sum\limits_{k=1}^{d} (p(\overline{f}f') + p(\overline{g}g')).e_k \otimes u_{ki}^2
		+ (q(\overline{f}f') + q(\overline{g}g')).e_k \otimes u_{ki} \\
	& = \alpha (\langle (f,g).e_i  | (f',g').e_j \rangle_A).
\end{align*}
Then by similar calculations as the ones done on $\alpha$, we obtain that $(\alpha,\beta)$ is a coaction of $\cA_h(d)$ on $E$.
Furthermore we have for all $n \in  \ZZ$ and all $i \in  \{1,\ldots, d\}$,
\[
	\beta(e_{ni}) =   \sum\limits_{k=1}^{d} e_{nk} \otimes u_{ki}^{r(k)}
\]
since $p(\sin (\pi n \dot ) ) = \dfrac{1+(-1)^{n+1}}{2} \sin(\pi n \dot)$
and $q(\sin (\pi n \dot)) = \dfrac{1+(-1)^n}{2} \sin (\pi n \dot)$.
Consequently, we have $\beta(V_n) \subset V_n \odot \cA_h(d)$.

We have moreover
\[
	\beta(\xi_0) =   \sum\limits_{i=1}^{d} \beta((0,1).e_i)
		=   \sum\limits_{i,k=1}^{d} (0,1).e_k \otimes u_{ki}^2 + (0,0).e_k \otimes u_{ki}
		=   \sum\limits_{k=1}^{d} (0,1).e_k \otimes 1 = \xi_0 \otimes 1,
\]
and it is clear that $(J \otimes *)\circ \beta = \beta \circ J$.
Lastly, we have $\tau( p(f).e_i) = \tau( f.e_i) = \displaystyle \int_0^1 f$ and $\tau( q(f).e_i) = 0$
since $\displaystyle \int_0^1 f(1-x)dx =\int_0^1 f(x)dx $ for all $f \in C([0,1])$.
Thus
\[
	(\tau \otimes id) \circ \alpha (f.e_i) =   \sum\limits_{k=1}^{d} \tau(p(f).e_k)u_{ki}^2 + \tau(q(f).e_k)u_{ki}
	= \displaystyle \left(\int_0^1 f \right). \left(  \sum\limits_{k=1}^{d} u_{ki}^2\right) = \int_0^1 f = \tau(f.e_i),
\]
which lets us conclude that $(\alpha,\beta)$ is a filtration preserving coaction of $\cA_h(d)$ on $E$, so that 
the quantum symmetry group of $E$ is the hyperoctahedral quantum group.

\vspace{1em}
\textbf{``Quantum isometry groups'' of quotients of $\boldsymbol{[0,1]\times  \{1,\ldots, d\}}$ }

Now we look at what might be the quantum isometry group of $d$ segments, all joined together in the point 0 and in the point 1.
The idea is to check the existence of a universal object in the full subcategory of $\cC(E,\tau,(V_i)_{i \in \ZZ}, J, \xi_0)$, consisting
of the Woronowicz $C^*$-algebras $Q$ whose coaction $(\alpha,\beta)$ on $E$ satisfies $\alpha(B) \subset B \otimes Q$,
where $B \subset A = C(I)^d$ is the $C^*$-algebra of continuous functions on the quotient of $[0,1] \times  \{1,\ldots, d\}$,
where for all $i,j \in  \{1,\ldots, d\}$, the points $(0,i)$ and $(0,j)$ are identified, and  the points $(1,i)$ and $(1,j)$ are identified.
That is to say:
 \[
	B = \{(f_1, \ldots, f_d) \in C(I)^d \tq \forall i,j \in  \{1,\ldots, d\}, f_i(0) = f_j(0) \text{ and } f_i(1) = f_j(1)\}.
\]

First notice that there is a filtration preserving coaction $(\alpha_0, \beta_0)$ of $\cA_s(d) \otimes C(\ZZ_2)$ on $E$
(where $\cA_s(d)$ is the quantum permutation group on $d$ points~\cite{wang_symgroups})
 characterized for $f \in C_0(I)$, $g \in C(I)$ and $i \in  \{1,\ldots, d\}$ by:
\begin{align*}
	\alpha_0 (  g.e_i )
		&=   \sum\limits_{k=1}^{d} p(g).e_k \otimes v_{ki} \otimes 1 + q(g).e_k \otimes v_{ki} \otimes z \\
	 \text{ and }
	\beta_0((f,g).e_i)&  =   \sum\limits_{k=1}^{d} (q(f),p(g)).e_k \otimes v_{ki} \otimes 1 + (p(f),q(g)).e_k \otimes v_{ki} \otimes z,
\end{align*}
where the $v_{ij}$'s are the canonical generators of $\cA_s(d)$ and $z \in C(\ZZ_2)$ is the function such that
$z(0) = 1$ and $z(1) =-1$.

Now let $(\alpha, \beta)$ be a  filtration preserving coaction of a Woronowicz $C^*$-algebra $Q$ on $E$,
and satisfying $\alpha(B) \subset B \otimes Q$.
Recall that there is a multiplicative and unitary matrix $w \in \cM_d(Q)$ such that for all $i,n$, 
$\beta(e_{ni}) = \sum\limits_{k=1}^d e_{nk} \otimes w_{ki}^{r(n)},$
$\alpha(\cos(\pi n \dot).e_i) = \sum\limits_{k=1}^d \cos(\pi n \dot).e_k \otimes w_{ki}^{r(n)}$, and satisfying:
\begin{itemize}
	\item $w_{ij}^* = w_{ij}$,
	\item  for $j \neq k$, $w_{ij} w_{ik}  = w_{ji}w_{ki} = 0$,
	\item $\sum\limits_{l = 1}^n w_{il}^2 =  \sum\limits_{l = 1}^n w_{li}^2 = 1$.
\end{itemize}
This leads to the existence of a unique morphism of Woronowicz $C^*$-algebras $\mu : \cA_s(d) \rightarrow Q$
such that $\mu(v_{ij}) = w_{ij}^2$.
Moreover, the condition $\alpha(B) \subset B \otimes Q$ means that
\[
	\forall a \in B, \forall \varepsilon \in \{0,1\}, \forall i,j \in  \{1,\ldots, d\}, \;
	(ev_{(\varepsilon,i)} \otimes id) \circ \alpha(a) =  (ev_{(\varepsilon,j)} \otimes id) \circ \alpha(a).
\]
In particular, 
$(ev_{(0,i)} \otimes id) \circ \alpha(\cos(\pi \dot).1_A)
	 =  (ev_{(0,j)} \otimes id) \circ \alpha(\cos(\pi \dot).1_A)$
for all $i,j \in \{1, \ldots, d\}$. But we have
\begin{align*}
\alpha(\cos(\pi \dot).1_A) = \alpha\left(\cos(\pi \dot).\left(  \sum\limits_{k=1}^{d} e_k\right)\right) 
	& =   \sum\limits_{k,l=1}^{d} \cos(\pi \dot) e_l \otimes w_{lk}
		=  \sum\limits_{l=1}^{d} \cos(\pi \dot) e_l \otimes \left(  \sum\limits_{k=1}^{d} w_{lk}\right)\!,
\end{align*}
thus $(ev_{(0,i)} \otimes id) \circ \alpha(\cos(\pi \dot).1_A) =   \sum\limits_{k=1}^{d} w_{ki}$.
Therefore we must have $\sum\limits_{k=1}^{d} w_{ki} =\sum\limits_{k=1}^{d} w_{kj}$ for all
$i,j \in  \{1,\ldots, d\} $. A simple calculation shows that $\omega = \sum\limits_{k=1}^{d} w_{ki} $
satisfies $\omega^2 =1$, $\omega^* = \omega$ and $\Delta(\omega) = \omega \otimes \omega$, so
that there exists a unique morphism of Woronowicz $C^*$-algebras
$\nu : C(\ZZ_2) \rightarrow Q$ sending $z$ to $\omega$ (recall that $C(\ZZ_2)$ is isomorphic to the universal Woronowicz
$C^*$-algebra generated by a unitary self-adjoint element). Since $\omega$ clearly commutes with the $w_{ij}^2$'s,
the morphism $\mu \otimes \nu : \cA_s(d) \otimes C(\ZZ_2) \rightarrow Q$ is well defined, and it is easy to see that it intertwines
the coactions. Thus $\cA_s(d) \otimes C(\ZZ_2)$ is a universal object in the subcategory of
$\cC(E, \tau, (V_i)_{i \in \ZZ}, J, \xi_0)$ we considered.

Note that the quantum isometry groups of other quotients of $[0,1] \times  \{1,\ldots, d\}$ can be computed in a similar fashion.
For example, the universal object in the category of Woronowicz $C^*$-algebras coacting on $E$ in a filtration preserving way
that additionally preserve the $C^*$-algebra $C = \{(f_1, \ldots , f_d\} \tq \forall i,j \in  \{1,\ldots, d\}, f_i(0) = f_j(0) = f_i(1) = f_j(1)\}$
(respectively  the $C^*$-algebra $D = \{(f_1, \ldots , f_d\} \tq \forall i,j \in  \{1,\ldots, d\}, f_i(0) = f_j(0)\}$)
is $\cA_h(d)$ (respectively $\cA_s(d)$).

\vspace{1em}
\textbf{Acknowledgments} $-$ The author is very grateful to Julien Bichon and Jean-Marie Lescure for their time and encouragement, and
to Georges Skandalis for his helpful suggestions and comments.

 \begin{small}
\bibliographystyle{abbrv}
\bibliography{QSG_HilbertMod.bib}
\end{small}
\end{document}